\documentclass[pdflatex,sn-mathphys-num]{sn-jnl}

\usepackage{graphicx}%
\usepackage{multirow}%
\usepackage{amsmath,amssymb,amsfonts}%
\usepackage{amsthm}%
\usepackage{mathrsfs}%
\usepackage[title]{appendix}%
\usepackage{xcolor}%
\usepackage{textcomp}%
\usepackage{manyfoot}%
\usepackage{booktabs}%
\usepackage{algorithm}%
\usepackage{algorithmicx}%
\usepackage{algpseudocode}%
\usepackage{listings}%
\usepackage[T1]{fontenc}
\usepackage[utf8]{inputenc}

\theoremstyle{thmstyleone}%
\theoremstyle{thmstyletwo}%

\theoremstyle{thmstylethree}%

\newtheorem{theo}{Theorem}
\newtheorem{lem}[theo]{Lemma}
\newtheorem{prop}[theo]{Proposition}
\newtheorem{cor}[theo]{Corollary}
\newtheorem{rem}{Remark}
\newtheorem{exam}{Example}
\newtheorem{definiti}{Definition}

\newenvironment{dem}[1][Proof]{\noindent \textbf{#1.} }{\ \rule{0.5em}{0.5em}}
\begin{document}

\title[Article Title]{Article Title}
\title[Relaxation in infinite convex programming]{Relaxation in infinite convex programming under
Slater-type regularity conditions}


\author[1]{\fnm{Rafael} \sur{Correa}}\email{rcorrea@dim.uchile.cl}
\equalcont{These authors contributed equally to this work.}

\author*[2]{\fnm{Abderrahim} \sur{Hantoute}}\email{hantoute@ua.es}
\equalcont{These authors contributed equally to this work.}

\author[2]{\fnm{Marco A.} \sur{López}}\email{marco.antonio@ua.es}
\equalcont{These authors contributed equally to this work.}

\affil[1]{\orgdiv{Mathematical Engineering Department}, \orgname{University of Chile}, \orgaddress{
\city{Santiago}, 
\country{Chile}}}

\affil*[2,3]{\orgdiv{Mathematical Department}, \orgname{University of Alicante}, \orgaddress{
\city{Alicante}, 
\country{Spain}}}



\abstract{The main purpose of this paper is to close the gap between the optimal values
of an infinite convex program and that of its biconjugate relaxation. It is
shown that Slater and continuity-type conditions guarantee such a zero-duality
gap. The approach uses calculus\ rules for the conjugation and
biconjugation\textbf{ }of\ the sum and pointwise supremum operations.\textbf{
}A second important objective of this work is to exploit these results on
relaxation by applying them in the context of duality theory.}

\keywords{Infinite optimization, convexity, biconjugate relaxation,
duality theory, conjugacy calculus, Slater conditions}


\pacs[MSC Classification (2010)]{26B05, 26J25, 49H05}

\maketitle

\section{Introduction}

The process of relaxation via the biconjugate function, also\ called
weak*-lower semicontinuous regularization, is widely used in robust
optimization, control theory, calculus of variations, equilibrium theory and
other topics (see, for instance, \cite{AuEk84}, \cite{BoVa88}, \cite{EkTe99},
\cite{HiLoVo11}, \cite{Ro70m}, etc.). In the context of variational calculus,
the use of the biconjugate relaxation allows for enlarging the feasible set
and reducing the optimal value of the original problem. Then, the existence of
relaxed (or weak) solutions contributes to have more information on the
solutions of the original solution.

Given an arbitrary family of proper convex functions $f_{t}:X\rightarrow
\mathbb{R}\cup\{+\infty\},$ $t\in T\cup\{0\}$, defined on a Banach space\ $X,$
we consider the convex optimization problem
\begin{equation}
(\mathcal{P})\text{ }\left\{
\begin{array}
[c]{l}%
\inf~~f_{0}(x)\\
\text{s.t. }f_{t}(x)\leq0,\text{ }t\in T,\\
x\in X.
\end{array}
\right.  \label{problemP}%
\end{equation}
When\ $T$ is infinite and $X$ is the Euclidean space, we are in the framework
of semi-infinite convex programming (see, e.g., \cite{LoSt07} and references
therein). The optimal value of $(\mathcal{P})$ is denoted by $v(\mathcal{P}%
)$,\ with the convention $v(\mathcal{P})=+\infty$ if $(\mathcal{P})$ is
infeasible (i.e., the feasible set $F(\mathcal{P})$ of $(\mathcal{P})$ is
empty).\ Associated with $(\mathcal{P}),$ we introduce its \emph{biconjugate
relaxation} posed in the bidual space $X^{\ast\ast}$
\begin{equation}
(\mathcal{P}^{\ast\ast})\text{ }\left\{
\begin{array}
[c]{l}%
\inf~~f_{0}^{\ast\ast}(z)\\
\text{s.t. }f_{t}^{\ast\ast}(z)\leq0,\text{ }t\in T,\\
z\in X^{\ast\ast},
\end{array}
\right.  \label{Problema2}%
\end{equation}
where $f_{t}^{\ast\ast}$ denotes\ the Fenchel biconjugate of $f_{t},$ $t\in
T.$

Since\ $f_{t}^{\ast\ast}\leq f_{t}$ for all $t\in T\cup\{0\},$ a weak
duality-like inequality always holds between $(\mathcal{P})$ and
$(\mathcal{P}^{\ast\ast}),$ that is
\begin{equation}
v(\mathcal{P}^{\ast\ast})\leq v(\mathcal{P}).\label{wi}%
\end{equation}
If $X$ is reflexive and, additionally, each $f_{t}$ is is lower
semicontinuous, then the Fenchel-Moreau-Rockafellar theorem ensures that
$f_{t}^{\ast\ast}=f_{t}$, for all $t\in T\cup\{0\},$ and thus $(\mathcal{P})$
and $(\mathcal{P}^{\ast\ast})$ coincide. Obviously, this problem is
meaningful\ only in the nonreflexive setting.

Inequality (\ref{wi}) however may be strict as we show through nontrivial
examples of $(\mathcal{P})$ provided in Section \ref{sec4}. Furthermore, it is
seen in Example \ref{examm} that, in any nonreflexive dual Banach space, there
exists a linear program with infinite constraints that exhibits a gap with its
biconjugate relaxation.

Our analysis is focused on guaranteeing that both, the original problem and an
adequate relaxation have the same optimal values. To this aim, we introduce a
reinforced alternative to $(\mathcal{P}^{\ast\ast});$ namely, $(\mathcal{P}%
_{\infty}^{\ast\ast}).$ We prove that, under Slater and continuity type
conditions, there is no gap between $(\mathcal{P})$ and $(\mathcal{P}_{\infty
}^{\ast\ast}).$ Other variants of $(\mathcal{P}_{\infty}^{\ast\ast})$ are
shown to be also useful for our purposes.

Equality
\begin{equation}
v(\mathcal{P}_{\infty}^{\ast\ast})=v(\mathcal{P}). \label{e1}%
\end{equation}
can also be viewed as a dual result and, therefore, it is related to Fenchel
duality in convex programming. In fact, the biconjugate relaxation
$(\mathcal{P}_{\infty}^{\ast\ast})$ provides a lower estimate for many
Fenchel-type duals of problem $(\mathcal{P}).\ $As a consequence of that, we
prove that conditions ensuring equality (\ref{e1}) also constitute sufficient
conditions for a zero-duality gap between $(\mathcal{P})$ and these Fenchel duals.

The paper is organized as follows. Notation and preliminary results are
gathered in Section \ref{sec1}. Section \ref{sec2} contains the main results
establishing\ zero-duality gap for several biconjugate relaxations. These
results are then applied in Section \ref{sec3} to optimization models
involving concave-like families of constraints. In Section \ref{sec4a}, we
explore the connections between biconjugate relaxation and Fenchel duality in
infinite convex programming. Section \ref{sec4} presents some illustrative
examples motivating\ the use of the different biconjugate relaxations and
showing the limitations of the standard one. Finally, concluding remarks are
given in Section \ref{sec5}.

\section{Notation and preliminary results\label{sec1}}

In the paper $X$ is a real Banach space, $X^{\ast}$ its dual, and $X^{\ast
\ast}$ its bidual. Given $x\in X$ and $x^{\ast}\in X^{\ast},$ we use the
notation $\langle x^{\ast},x\rangle\equiv\left\langle x,x^{\ast}\right\rangle
:=x^{\ast}(x).$ We consider the inclusion $X\subset X^{\ast\ast}$ by
identifying $X\ni x\leftrightarrow\left\langle \cdot,x\right\rangle \in
X^{\ast\ast}.$ The zero vector is denoted by $\theta,$ and $B_{X}$ stands for
the closed unit ball of $X.$ We use $\left\Vert \cdot\right\Vert $ to denote
the norm in any Banach space, in absence of any possible\textbf{ }confusion.
For instance, $\ell_{1}$ is the Banach space of real sequences $(x_{n})_{n}$
such that $\left\Vert (x_{n})_{n}\right\Vert :=\sum_{n\geq1}\left\vert
x_{n}\right\vert <+\infty,$ while\ $\ell_{\infty}$ is the dual of $\ell_{1},$
endowed with the norm $\left\Vert (x_{n})_{n}\right\Vert :=\sup_{n\geq
1}\left\vert x_{n}\right\vert $. By $\ell_{1}^{+}$ and $\ell_{\infty}^{+}$ we
represent the sets of nonnegative sequences in $\ell_{1}$ and $\ell_{\infty},$
respectively. We denote\ $\overline{\mathbb{R}}:=\mathbb{R}\cup\{-\infty
,+\infty\}$ and $\mathbb{R}_{\infty}:=\mathbb{R}\cup\{+\infty\}$, with the
convention $0\cdot(+\infty)=+\infty.$

Given a set $A\subset X,$ by $\operatorname*{co}(A),$ $\operatorname*{int}%
(A),$ and $\operatorname*{cl}(A)$ (or\ $\overline{A})$ we denote the convex
hull, the interior, and the closure of $A$, respectively. We set
$\overline{\operatorname*{co}}(A):=\operatorname*{cl}(\operatorname*{co}(A)).$

Given a function $f:X\longrightarrow\overline{\mathbb{R}}$, the sets
$\operatorname*{dom}f:=\{x\in X:\ f(x)<+\infty\}$ and $\operatorname*{epi}%
f:=\{(x,\lambda)\in X\times\mathbb{R}:\ f(x)\leq\lambda\}$ are\ the
(effective) domain and epigraph of $f$, respectively. We set $[f\leq
\alpha]:=\{x\in X:f(x)\leq\alpha\},$ $\alpha\in\mathbb{R}$. The positive part
of $f$ is the function $f^{+}:=\max\{f,0\}.$ In particular, the indicator
function of $A$ is defined by $\mathrm{I}_{A}(x)=0$ if $x\in A,$ and
$\mathrm{I}_{A}(x)=+\infty$ otherwise. The function\ $f$ is proper if
$f\not \equiv +\infty$ and $f>-\infty$, lower semicontinuous (lsc, for short)
if $\operatorname*{epi}f$ is closed, and convex if $\operatorname*{epi}f$ is
convex.\ The closed and the closed convex hulls of $f$ are, respectively,\ the
functions $\operatorname*{cl}f$ (or $\bar{f})$ and $\overline
{\operatorname*{co}}f$ such that\ $\operatorname*{epi}\bar{f}%
=\operatorname*{cl}(\operatorname*{epi}f)$ and $\operatorname*{epi}%
(\overline{\operatorname*{co}}f)=\overline{\operatorname*{co}}%
(\operatorname*{epi}f).$ The family of proper, convex, and lsc functions
defined on $X$ is denoted by $\Gamma_{0}(X).$

The conjugate of $f$ is the lsc convex function $f^{\ast}:X^{\ast}%
\rightarrow\overline{\mathbb{R}}\ $defined by
\[
f^{\ast}(x^{\ast}):=\sup\{\left\langle x^{\ast},x\right\rangle -f(x),\ x\in
X\},
\]
while the biconjugate of $f$ is the conjugate of $f^{\ast};$ that is,
\[
f^{\ast\ast}(z):=\sup\{\left\langle z,x^{\ast}\right\rangle -f^{\ast}(x^{\ast
}),\ x^{\ast}\in X^{\ast}\},\text{ }z\in X^{\ast\ast}.
\]
The support function of a nonempty set $A\subset X$ is $\mathrm{\sigma}%
_{A}:=(\mathrm{I}_{A})^{\ast}.$ Due to the Fenchel-Moreau-Rockafellar theorem
(\cite[Theorem 3.2.2(ii)]{CHLBook}), provided that the convex function $f$
admits a continuous affine minorant (for instance, when $f\in\Gamma_{0}(X)$)
we have
\begin{equation}
f^{\ast\ast}(z)=\overline{f}^{w^{\ast\ast}}(z):=\liminf_{x\rightharpoonup
z,\text{ }x\in X}f(x),\text{ for all }z\in X^{\ast\ast}, \label{moreau}%
\end{equation}
where \textquotedblleft$\rightharpoonup$\textquotedblright\ stands for the
convergence in\ the weak-* topology of\ $X^{\ast\ast},$ which is denoted by
$w^{\ast\ast}$ (usually represented by $\sigma(X^{\ast\ast},X^{\ast})$), and
\[
\liminf_{x\rightharpoonup z,\text{ }x\in X}f(x):=\sup_{x^{\ast}\in X^{\ast
},\text{ }\varepsilon>0\text{ }}\inf_{x\in B_{x^{\ast}}(z,\varepsilon)\text{
}}f(x),
\]
where $B_{x^{\ast}}(z,\varepsilon):=\{x\in X:\left\vert \left\langle x^{\ast
},x\right\rangle -\left\langle z,x^{\ast}\right\rangle \right\vert
\leq\varepsilon\}.$ Sometimes we write $\operatorname*{cl}\nolimits^{w^{\ast
\ast}}(f)$ instead of $\overline{f}^{w^{\ast\ast}}.$ When $f$ is lsc, in $X$
we have $f=f^{\ast\ast}.$

For every nonempty\ closed set $A\subset X$ $(\subset X^{\ast\ast})$ we have
\[
\mathrm{I}_{A}^{\ast\ast}=\operatorname*{cl}\nolimits^{w^{\ast\ast}%
}(\mathrm{I}_{A})=\mathrm{I}_{\overline{A}^{w^{\ast\ast}}},
\]
where $\overline{A}^{w^{\ast\ast}}$ denotes the $w^{\ast\ast}$-closure of $A$
in $X^{\ast\ast}.$ Thus, as a consequence, we obtain
\begin{equation}
\overline{\operatorname*{dom}f^{\ast\ast}}^{w^{\ast\ast}}=\overline
{\operatorname*{dom}f}^{w^{\ast\ast}}. \label{sets}%
\end{equation}

Given the\textbf{ }convex functions $f_{0},f_{1},$ $\cdots,$ $f_{m}%
:X\rightarrow\mathbb{R}_{\infty}$ such that $f_{1},$ $\cdots,$ $f_{m}$ are
continuous at some point in $\operatorname*{dom}f_{0}$, we have (\cite{Mo66},
\cite{Ro70m})\
\begin{equation}
(f_{0}+f_{1}+\cdots+f_{m})^{\ast\ast}=f_{0}^{\ast\ast}+f_{1}^{\ast\ast}%
+\cdots+f_{m}^{\ast\ast}.\label{sumbicon}%
\end{equation}
The same hypothesis guarantees that
\begin{equation}
(\max\{f_{0},f_{1},\cdots,f_{m}\})^{\ast\ast}=\max\{f_{0}^{\ast\ast}%
,f_{1}^{\ast\ast},\cdots,f_{m}^{\ast\ast}\}.\label{maxbicon}%
\end{equation}
The last relation also holds under the Brézis-Attouch condition (\cite{FS00}).
Other conditions have been proposed in \cite{BW08, CHL26RACSAM, Za08}. In
particular, combining (\ref{maxbicon}) and (\ref{sets}), the continuity
assumption above implies that (\cite[Corollary 7]{CHL26-JOTA})
\begin{equation}
\operatorname*{cl}\nolimits^{w^{\ast\ast}}(\operatorname*{dom}(f^{\ast\ast
}))=\operatorname*{cl}\nolimits^{w^{\ast\ast}}(\operatorname*{dom}%
f)=\cap_{k=1}^{m}\operatorname*{cl}\nolimits^{w^{\ast\ast}}%
(\operatorname*{dom}f_{k})=\cap_{k=1}^{m}\operatorname*{cl}\nolimits^{w^{\ast
\ast}}(\operatorname*{dom}f_{k}^{\ast\ast}).\label{hamdo}%
\end{equation}

Finally, let\ $f_{0},f:X\rightarrow\mathbb{R}_{\infty}$ be two convex
functions such that $\mu:=\inf_{f(x)\leq0}f_{0}(x)\in\mathbb{R}$, and
$f(x_{0})<0$ for some $x_{0}\in\operatorname*{dom}f_{0}$ (i.e., Slater
condition holds).$\ $Then, a real number $\alpha\in\mathbb{R}$ satisfies
\begin{equation}
0=\inf_{x\in X}\max\{f_{0}(x)-\alpha,\text{ }f(x)\}\label{ca}%
\end{equation}
if and only if $\alpha=\mu.$ We refer, e.g., to \cite[Lemma 3.7]{Ca23} for a
proof in finite-dimensional setting; the argument extends readily to our framework.

\section{Relaxation theory\label{sec2}}

In this section, we present the main results of this work dealing with the
relaxation of the convex optimization problem $(\mathcal{P})$ defined in
(\ref{problemP}). We provide several results relating the optimal value of
$(\mathcal{P})$ to the optimal values\textbf{ }of various relaxations.

We suppose that $T$\ is countable, say $T:=\{1,2,\cdots\}$, in which case the
problem $(\mathcal{P})$ takes the form\
\[
(\mathcal{P})\text{ \ }\left\{
\begin{array}
[c]{l}%
\inf~~f_{0}(x)\\
\text{s.t. }f_{k}(x)\leq0,\text{ }k\geq1,\\
\text{ \ \ \ \ }x\in X,
\end{array}
\right.
\]
where $f_{k}:X\rightarrow\mathbb{R}_{\infty},$ $k\geq0$, are assumed to be
proper and convex. We further assume that the optimal value $v(\mathcal{P}%
)<+\infty;$ equivalently, the feasible set $F(\mathcal{P})$ is nonempty.$\ $

\begin{rem}
If $T$ is not countable but $X$ is separable, then $v(\mathcal{P})$ remains
unchanged if we replace $T$ by an appropriate countable subset (see,
e.g.,\ \cite[Lemma 3.1]{HaKrLo25b}).\ 
\end{rem}

As one can expect, the following Slater condition will be crucial in our analysis.

\begin{definiti}
\label{DefSlater} We say that the system\ of constraints\textbf{ }%
$\{f_{t}(x)\leq0,$ $t\in T\}$ satisfies the\ (strong) Slater condition, if
there exists $x_{0}\in X$ such that
\[
\sup_{t\in T}f_{t}(x_{0})<0.
\]
The point\ $x_{0}$ is called a Slater point of the system $\{f_{t}(x)\leq0,$
$t\in T\}.$ When, additionally, $x_{0}\in\operatorname*{dom}f_{0},$ we say
that problem $(\mathcal{P})$ satisfies the Slater condition.
\end{definiti}

A natural\ candidate for the relaxation of $(\mathcal{P})$ is given by
\begin{equation}
(\mathcal{P}^{\ast\ast})\text{ \ }%
\begin{array}
[c]{l}%
\inf~~f_{0}^{\ast\ast}(z)\\
\text{s.t. }f_{k}^{\ast\ast}(z)\leq0,\text{ }k\geq1,\\
\text{ \ \ \ \ }z\in X^{\ast\ast}.
\end{array}
\label{eti}%
\end{equation}
It is straightforward\ that $v(\mathcal{P}^{\ast\ast})\leq v(\mathcal{P}).$
Moreover, by the Fenchel-Moreau-Rockafellar theorem (\ref{moreau}), when\ $X$
is reflexive and the $f_{k}$'s are lsc, we obviously have zero gap,
$v(\mathcal{P}^{\ast\ast})=v(\mathcal{P}),$ since both problems coincide.

The following proposition shows that a zero gap between $(\mathcal{P})$ and
$(\mathcal{P}^{\ast\ast})$ also holds when the number of constraints in
$(\mathcal{P})$ is finite.\ 

\begin{prop}
\label{prop2}Assume that the convex functions $f_{k},$ $k=1,\cdots,m$, are
finite and continuous at some point in $\operatorname*{dom}f_{0}.$ Then, under
the Slater condition, we have\
\[
v(\mathcal{P})=v(\mathcal{P}^{\ast\ast}).
\]

\end{prop}

\begin{dem}
We may assume that $v(\mathcal{P})\in\mathbb{R};$ otherwise, we are done
because the relations $f_{0}^{\ast\ast}\leq f_{0}$ and $f_{k}^{\ast\ast}\leq
f_{k}^{\ast\ast},$ for $k\geq1,$ imply
\[
v(\mathcal{P}^{\ast\ast})\leq\inf_{f_{k}^{\ast\ast}(x)\leq0,\text{ }%
k\geq1,\text{ }x\in X\text{\ }}f_{0}^{\ast\ast}(x)\leq\inf_{f_{k}%
(x)\leq0,\text{ }k\geq1,\text{ }x\in X\text{ \ }}f_{0}(x)=v(\mathcal{P}).
\]
Define the performance function
\[
g:=\max\{f_{0}-v(\mathcal{P}),\ f_{k},\ k=1,\cdots,m\}.
\]
Then, by\ \cite[Proposition 5.2.4(ii)]{CHLBook},
\[
\operatorname*{cl}g:=\max\{(\operatorname*{cl}f_{0})-v(\mathcal{P}%
),\ \operatorname*{cl}f_{k},\ k=1,\cdots,m\}.
\]
Using (\ref{maxbicon}),\ the current assumption implies that
\begin{align*}
g^{\ast\ast}=(\operatorname*{cl}g)^{\ast\ast} &  =\max\{(\operatorname*{cl}%
f_{0})^{\ast\ast}-v(\mathcal{P}),\ (\operatorname*{cl}f_{k})^{\ast\ast
},\ k=1,\cdots,m\}\\
&  =\max\{f_{0}^{\ast\ast}-v(\mathcal{P}),\ f_{k}^{\ast\ast},\ k=1,\cdots,m\},
\end{align*}
and, so,
\begin{align}
0=\inf_{X}g(x)=\inf_{X^{\ast\ast}}g^{\ast\ast}(z) &  =\inf_{X^{\ast\ast}}%
\max\{f_{0}^{\ast\ast}(z)-v(\mathcal{P}),\ f_{k}^{\ast\ast}(z),\ k=1,\cdots
,m\}\nonumber\\
&  =\inf_{X^{\ast\ast}}\max\{f_{0}^{\ast\ast}(z)-v(\mathcal{P}),\ \max_{1\leq
k\leq m}f_{k}^{\ast\ast}(z)\}.\label{mm}%
\end{align}
By the Slater assumption, there exists $x_{0}\in\operatorname*{dom}%
f_{0}\subset\operatorname*{dom}f_{0}^{\ast\ast}$ such that
\[
\max_{k\geq1}f_{k}^{\ast\ast}(x_{0})\leq\max_{k\geq1}f_{k}(x_{0})<0.
\]
Therefore, applying\ (\ref{ca}) to the problem $\inf_{\max\limits_{1\leq k\leq
m}f_{k}^{\ast\ast}(z)\leq0,\text{ }z\in X^{\ast\ast}\text{  }}f_{0}^{\ast\ast
}(z)$, (\ref{mm}) entails that $v(\mathcal{P})=v(\mathcal{P}^{\ast\ast}).$
\end{dem}

As evidenced in Example \ref{examm} in Section \ref{sec4}, when problem
$(\mathcal{P})$ involves infinitely many constraints, the biconjugate
relaxation $(\mathcal{P}^{\ast\ast})$ introduced above needs to be reinforced
to eliminate the gap with $(\mathcal{P})$. To overcome this issue, we
associate with $(\mathcal{P})$ a strengthened relaxed problem,
\begin{equation}
(\mathcal{P}_{\infty}^{\ast\ast})\text{ }\left\{
\begin{array}
[c]{l}%
\inf~~f_{0}^{\ast\ast}(z)\\
\text{s.t. }f_{k}^{\ast\ast}(z)\leq0,\text{ }k\geq1,\text{ }\\
\text{ \ \ \ \ }f_{\infty}^{\ast\ast}(z)\leq0,\\
\text{ \ \ \ \ }z\in\overline{\operatorname*{dom}(f^{\ast\ast})}^{w^{\ast\ast
}},
\end{array}
\right.  \label{pdual}%
\end{equation}
where $f:=\sup_{k\geq1}f_{k},$ the function $f_{\infty}:X\rightarrow
\mathbb{R}_{\infty}$ is defined as
\begin{equation}
f_{\infty}(x):=\limsup_{k\rightarrow\infty}f_{k}(x),\text{ } \label{fs}%
\end{equation}
and $f_{\infty}^{\ast\ast}$ is its biconjugate. When the number of constraints
in $(\mathcal{P})$ is finite, we set $f_{\infty}\equiv-\infty$.

Next we show that, in the case of convex optimization problems with finitely
many constraints, this formulation reduces to the relaxation $(\mathcal{P}%
^{\ast\ast})$ introduced in (\ref{eti}).

\begin{prop}
If\ the number of constraints in $(\mathcal{P})$ is finite, say $m,$ and the
functions $f_{k},$ $k=1,\cdots,m$, are finite and continuous at some point in
$\operatorname*{dom}f_{0}$, then $(\mathcal{P}_{\infty}^{\ast\ast})$ coincides
with\ $(\mathcal{P}^{\ast\ast})$.
\end{prop}

\begin{dem}
On the one hand, we have $f_{\infty}\equiv-\infty,$ so that $f_{\infty}%
^{\ast\ast}\equiv-\infty$ and the inequality $f_{\infty}^{\ast\ast}(z)\leq0$
in $(\mathcal{P}_{\infty}^{\ast\ast})$ trivially holds. On the\ other hand,
the continuity assumption implies that $f^{\ast\ast}=\operatorname*{cl}%
\nolimits^{w^{\ast\ast}}(f),$ by (\ref{moreau}), and
\[
\operatorname*{cl}\nolimits^{w^{\ast\ast}}(\operatorname*{dom}(f^{\ast\ast
}))=\operatorname*{cl}\nolimits^{w^{\ast\ast}}(\operatorname*{dom}%
f)=\operatorname*{cl}\nolimits^{w^{\ast\ast}}(\cap_{k=1}^{m}%
\operatorname*{dom}f_{k}).
\]
Thus, by (\ref{hamdo})\textbf{\ }we obtain
\begin{equation}
\operatorname*{cl}\nolimits^{w^{\ast\ast}}(\operatorname*{dom}(f^{\ast\ast
}))=\cap_{k=1}^{m}\operatorname*{cl}\nolimits^{w^{\ast\ast}}%
(\operatorname*{dom}f_{k})=\cap_{k=1}^{m}\operatorname*{cl}\nolimits^{w^{\ast
\ast}}(\operatorname*{dom}f_{k}^{\ast\ast}), \label{ar}%
\end{equation}
and the constraint $z\in\operatorname*{cl}\nolimits^{w^{\ast\ast}%
}(\operatorname*{dom}(f^{\ast\ast}))$ in $(\mathcal{P}_{\infty}^{\ast\ast})$
is redundant. This shows that $(\mathcal{P}_{\infty}^{\ast\ast})$ coincides
with problem $(\mathcal{P}^{\ast\ast})$.\smallskip
\end{dem}

As expected, the optimal value $v(\mathcal{P}_{\infty}^{\ast\ast})$ does not
exceed $v(\mathcal{P})$; this is proved in the following lemma.

\begin{lem}
\label{wd}We have
\[
v(\mathcal{P}_{\infty}^{\ast\ast})\leq v(\mathcal{P}).
\]

\end{lem}

\begin{dem}
Since $X\subset X^{\ast\ast}$ and $f_{k}^{\ast\ast}(x)\leq f_{k}(x)$ for all
$x\in X$ and $k\geq0,$ it follows that\
\begin{align*}
v(\mathcal{P}_{\infty}^{\ast\ast}) &  \leq\inf_{\substack{f_{k}^{\ast\ast
}(x)\leq0,\text{ }k\geq1,\text{ }f_{\infty}^{\ast\ast}(x)\leq0\text{ }%
\\x\in\overline{\operatorname*{dom}f^{\ast\ast}}^{w^{\ast\ast}},\text{ }x\in
X\text{ \ \ \ \ \ \ \ \ }}}f_{0}^{\ast\ast}(x)\medskip\\
&  \leq\inf_{\substack{f_{k}(x)\leq0,\text{ }k\geq1,\text{ }f_{\infty}%
^{\ast\ast}(x)\leq0\text{ }\\x\in\overline{\operatorname*{dom}f^{\ast\ast}%
}^{w^{\ast\ast}},\text{ }x\in X\text{ \ \ \ \ \ }}}f_{0}(x)=\inf_{f_{k}%
(x)\leq0,\text{ }k\geq1,\text{ }x\in X}f_{0}(x)=v(\mathcal{P}).
\end{align*}
The first equality above comes from the fact that, if $x$ is feasible for
$(\mathcal{P}),$ then $x\in\operatorname*{dom}f\subset\operatorname*{dom}%
f^{\ast\ast}\subset\overline{\operatorname*{dom}f^{\ast\ast}}^{w^{\ast\ast}}$
and $f_{\infty}^{\ast\ast}(x)\leq f_{\infty}(x)=\limsup_{k\rightarrow\infty
}f_{k}(x)\leq0.\smallskip$
\end{dem}

We proceed\ by establishing the equality\ between the optimal values of
$(\mathcal{P})$ and its relaxed formulation $(\mathcal{P}_{\infty}^{\ast\ast
}),$ under the Slater condition and the continuity of the supremum function
$\sup_{k\geq1}f_{k}.$ The following lemma involves the notion of upper\ sum of
a family of functions\textbf{ }$f_{k}:X\rightarrow\mathbb{R}_{\infty},$
$k\geq1:$
\begin{equation}
\left(
{\textstyle\sum_{k\geq1}}
f_{k}\right)  (x):=\limsup_{m\rightarrow+\infty}%
{\textstyle\sum_{k=1,\cdots,m}}
f_{k}(x),\text{ }x\in X. \label{uppers}%
\end{equation}

\begin{lem}
\label{lemHKL} Suppose that $(\mathcal{P})$ satisfies the Slater condition and
that $v(\mathcal{P})\in\mathbb{R}$. Then, there exist $\lambda:=(\lambda
_{1},\lambda_{2},\cdots)\in\ell_{1}^{+},$ $\hat{\lambda}:=(\hat{\lambda}%
_{1},\hat{\lambda}_{2},\cdots)\in\ell_{\infty}^{+},$ and $\lambda_{\infty}%
\geq0$ such that
\begin{align}
v(\mathcal{P})  &  =\inf_{x\in X\text{ }}(f_{0}(x)+%
{\textstyle\sum_{k\geq1}}
\lambda_{k}f_{k}^{+}(x)+\lambda_{\infty}f_{\infty}(x))\label{a1b}\\
&  =\inf_{x\in X\text{ }}(f_{0}(x)+%
{\textstyle\sum_{k\geq1}}
\hat{\lambda}_{k}f_{k}^{+}(x)). \label{a2}%
\end{align}

\end{lem}

\begin{dem}
Relation (\ref{a2}) is established in \cite[Lemma 3.1(ii)]{HaKrLo25b}.
Moreover, by \cite[Theorem 3.1]{HaKrLo25b}, the Slater condition gives rise to
the existence of\ some $\lambda\in\ell_{1}^{+}$ and $\lambda_{\infty}\geq0$
such that
\begin{equation}
v(\mathcal{P})=\inf_{X\text{ }}(f_{0}+%
{\textstyle\sum_{k\geq1}}
\lambda_{k}f_{k}+\lambda_{\infty}f_{\infty}), \label{do}%
\end{equation}
where $%
{\textstyle\sum_{k\geq1}}
\lambda_{k}f_{k}$ is the upper sum defined in (\ref{uppers}). Therefore, since
$%
{\textstyle\sum_{k\geq1}}
\lambda_{k}f_{k}\leq%
{\textstyle\sum_{k\geq1}}
\lambda_{k}f_{k}^{+},$ we obtain
\begin{align*}
v(\mathcal{P})  &  \leq\inf_{X\text{ }}(f_{0}+%
{\textstyle\sum_{k\geq1}}
\lambda_{k}f_{k}^{+}+\lambda_{\infty}f_{\infty})\\
&  \leq\inf_{f_{k}(x)\leq0,\text{ }k\geq1\text{ }}(f_{0}(x)+\lambda_{\infty
}f_{\infty}(x))\leq\inf_{f_{k}(x)\leq0,\text{ }k\geq1\text{ }}f_{0}%
(x)=v(\mathcal{P}),
\end{align*}
which yields (\ref{a1b}).
\end{dem}

We also need the following lemma, which\ extends the biconjugate sum rule in
(\ref{sumbicon}) to upper\ sums.

\begin{lem}
\label{prop3} Consider a countable family\ of nonnegative convex functions
$f_{k}:X\rightarrow\mathbb{R}_{\infty}$ and let $\alpha_{k}>0,$ $k\geq1,$ such
that $(\alpha_{k})_{k}\in\ell_{1}$.$\ $Denote $f:=\sup_{k\geq1}f_{k}$ and
assume that\ each $f_{k}$ is continuous at some point in $\operatorname*{dom}%
f.$ Then\
\begin{equation}
\left(
{\textstyle\sum_{k\geq1}}
\alpha_{k}f_{k}\right)  ^{\ast\ast}(z)=%
{\textstyle\sum_{k\geq1}}
\alpha_{k}f_{k}^{\ast\ast}(z),\text{ for all }z\in\operatorname*{dom}%
f^{\ast\ast}. \label{rr}%
\end{equation}

\end{lem}

\begin{dem}
Denote $g:=%
{\textstyle\sum_{k\geq1}}
\alpha_{k}f_{k}.$ Since $f_{k}\geq f_{k}^{\ast\ast},$ for each $k\geq1,$ we
have $g\geq%
{\textstyle\sum_{k\geq1}}
\alpha_{k}f_{k}^{\ast\ast}.$ Because the functions $f_{k}^{\ast\ast}$ are also
nonnegative, we have
\[%
{\textstyle\sum_{k\geq1}}
\alpha_{k}f_{k}^{\ast\ast}=\sup_{m\geq1}%
{\textstyle\sum_{k=1,\cdots,m}}
\alpha_{k}f_{k}^{\ast\ast},
\]
and\ the function $\sum_{k\geq1}\alpha_{k}f_{k}^{\ast\ast}$ is convex and lsc.
Moreover, for every $x\in\operatorname*{dom}f,$ we have
\[
\sum_{k\geq1}\alpha_{k}f_{k}^{\ast\ast}(x)\leq\sum_{k\geq1}\alpha_{k}%
f_{k}(x)\leq\sum_{k\geq1}\alpha_{k}f(x)<+\infty,
\]
so that $\sum_{k\geq1}\alpha_{k}f_{k}^{\ast\ast}\in\Gamma_{0}(X^{\ast\ast})$,
entailing
\[
g^{\ast\ast}\geq(%
{\textstyle\sum_{k\geq1}}
\alpha_{k}f_{k}^{\ast\ast})^{\ast\ast}=%
{\textstyle\sum_{k\geq1}}
\alpha_{k}f_{k}^{\ast\ast},
\]
i.e., the inequality \textquotedblleft$\geq$\textquotedblright\ in (\ref{rr}).

To prove the converse inequality, we fix\ $\varepsilon>0$ and choose
$k_{0}\geq1$ large enough to satisfy
\[%
{\textstyle\sum_{k>m}}
\alpha_{k}\leq%
{\textstyle\sum_{k>k_{0}}}
\alpha_{k}\leq\varepsilon,\text{ for all }m\geq k_{0}.
\]
Fix\ $m\geq k_{0}.$ We have
\begin{equation}
g:=\limsup_{n\rightarrow\infty}%
{\textstyle\sum_{1\leq k\leq n}}
\alpha_{k}f_{k}=%
{\textstyle\sum_{1\leq k\leq m}}
\alpha_{k}f_{k}+%
{\textstyle\sum_{k>m}}
\alpha_{k}f_{k}\leq%
{\textstyle\sum_{1\leq k\leq m}}
\alpha_{k}f_{k}+\varepsilon f.\label{he}%
\end{equation}
For each $k\geq1,$ let\ $x_{k}\in\operatorname*{dom}f$ $(\subset\cap_{j\geq
1}\operatorname*{dom}f_{j})$ be\ a continuity point of $f_{k}.$ Then
\[
\bar{x}=%
{\textstyle\sum_{1\leq k\leq m}}
\frac{x_{k}}{m}\in\cap_{j=1,\cdots,m}\operatorname*{int}(\operatorname*{dom}%
f_{j}),\
\]
and the function $%
{\textstyle\sum_{1\leq k\leq m}}
\alpha_{k}f_{k}$ is continuous at $\bar{x}\in\operatorname*{dom}f.$ Hence,
(\ref{he}) together with (\ref{sumbicon}) yields
\[
g^{\ast\ast}(z)\leq%
{\textstyle\sum_{1\leq k\leq m}}
\alpha_{k}f_{k}^{\ast\ast}(z)+\varepsilon f^{\ast\ast}(z),\text{ for all
}m\geq k_{0}.
\]
Finally, the conclusion follows when $\varepsilon$ goes to zero and
$z\in\operatorname*{dom}f^{\ast\ast}.\smallskip$
\end{dem}

We are now ready to establish\ the main result of this section.\ 

\begin{theo}
\label{main1}$\ $Suppose that $(\mathcal{P})$ satisfies the Slater condition.
If the function $f:=\sup_{k\geq1}f_{k}$ is continuous at some point in
$\operatorname*{dom}f_{0},$ then
\[
v(\mathcal{P}_{\infty}^{\ast\ast})=v(\mathcal{P}).
\]

\end{theo}

\begin{dem}
We may assume that $v(\mathcal{P})\in\mathbb{R};$ otherwise, $v(\mathcal{P}%
)=-\infty$ and we are done because $v(\mathcal{P}_{\infty}^{\ast\ast})\leq
v(\mathcal{P})$ by Lemma \ref{wd}.$\ $Moreover, by the same Lemma \ref{wd}, we
only need to prove the inequality $v(\mathcal{P})\leq v(\mathcal{P}_{\infty
}^{\ast\ast})$. According to Lemma \ref{lemHKL}, the Slater condition yields
some\ $\lambda:=(\lambda_{1},\lambda_{2},\cdots)\in\ell_{1}^{+}$ and
$\lambda_{\infty}\geq0$ such that
\[
v(\mathcal{P})=\inf_{X\text{ }}(f_{0}+%
{\textstyle\sum_{k\geq1}}
\lambda_{k}f_{k}^{+}+\lambda_{\infty}f_{\infty}).
\]
Equivalently, since a function and its biconjugate have the same infimum, we
obtain
\begin{equation}
v(\mathcal{P})=\inf_{X^{\ast\ast}\text{ }}(f_{0}+%
{\textstyle\sum_{k\geq1}}
\lambda_{k}f_{k}^{+}+\lambda_{\infty}f_{\infty})^{\ast\ast}. \label{zay}%
\end{equation}
Thus, the proof reduces to decomposing the\ biconjugate of the sums appearing
in the last equation.

Let\ $\bar{x}\in\operatorname*{dom}f\cap\operatorname*{dom}f_{0}$ be a
continuity point of $f,$ and denote $g:=%
{\textstyle\sum_{k\geq1}}
\lambda_{k}f_{k}^{+}.$ Since $g\leq(%
{\textstyle\sum_{k\geq1}}
\lambda_{k})f^{+}$ and $f_{\infty}\leq f,$ both $g$ and $f_{\infty}$ are
continuous at $\bar{x}\in\operatorname*{dom}f_{0}.$ Hence, $\lambda_{\infty
}f_{\infty}$ is clearly continuous at $\bar{x}$ when $\lambda_{\infty}>0.$ In
addition, when $\lambda_{\infty}=0,$ our convention $0\cdot(+\infty)=+\infty$
entails $\lambda_{\infty}f_{\infty}=\mathrm{I}_{\operatorname*{dom}f_{\infty}%
}\leq\mathrm{I}_{\operatorname*{dom}f},$ and therefore $\lambda_{\infty
}f_{\infty}$ is also continuous at $\bar{x}.$ Thus, using (\ref{sumbicon}),
(\ref{zay}) implies
\begin{equation}
v(\mathcal{P})=\inf_{X^{\ast\ast}\text{ }}(f_{0}^{\ast\ast}+g^{\ast\ast
}+(\lambda_{\infty}f_{\infty})^{\ast\ast}).\label{abo}%
\end{equation}
Observe that $(\lambda_{\infty}f_{\infty})^{\ast\ast}=\lambda_{\infty
}f_{\infty}^{\ast\ast}$ when $\lambda_{\infty}>0.$ Moreover, when
$\lambda_{\infty}=0$ we have
\[
(\lambda_{\infty}f_{\infty})^{\ast\ast}=(\mathrm{I}_{\operatorname*{dom}%
f_{\infty}})^{\ast\ast}=\mathrm{I}_{\overline{\operatorname*{dom}f_{\infty}%
}^{w^{\ast\ast}}}\leq\mathrm{I}_{\operatorname*{dom}f_{\infty}^{\ast\ast}%
}=0\cdot f_{\infty}^{\ast\ast}.
\]
So, in all cases we have $(\lambda_{\infty}f_{\infty})^{\ast\ast}\leq
\lambda_{\infty}f_{\infty}^{\ast\ast}$ and (\ref{abo})\ yields
\begin{equation}
v(\mathcal{P})\leq\inf_{X^{\ast\ast}\text{ }}(f_{0}^{\ast\ast}+g^{\ast\ast
}+\lambda_{\infty}f_{\infty}^{\ast\ast}).\label{re0}%
\end{equation}
To simplify $g^{\ast\ast}$ we observe that, since $\operatorname*{dom}%
f\subset\operatorname*{dom}f_{k}$ for all $k\geq1,$\
\[
g=%
{\textstyle\sum_{k\geq1,\text{ }\lambda_{k}>0}}
\lambda_{k}f_{k}^{+}+%
{\textstyle\sum_{k\geq1,\text{ }\lambda_{k}=0}}
\mathrm{I}_{\operatorname*{dom}f_{k}}\leq%
{\textstyle\sum_{k\geq1,\text{ }\lambda_{k}>0}}
\lambda_{k}f_{k}^{+}+\mathrm{I}_{\operatorname*{dom}f}.
\]
Applying again (\ref{sumbicon}), we obtain
\begin{equation}
g^{\ast\ast}\leq(%
{\textstyle\sum_{k\geq1,\text{ }\lambda_{k}>0}}
\lambda_{k}f_{k}^{+}+\mathrm{I}_{\operatorname*{dom}f})^{\ast\ast}=(%
{\textstyle\sum_{k\geq1,\text{ }\lambda_{k}>0}}
\lambda_{k}f_{k}^{+})^{\ast\ast}+\mathrm{I}_{\overline{\operatorname*{dom}%
f}^{w^{\ast\ast}}}.\label{re}%
\end{equation}
Furthermore, by Lemma \ref{prop3} and the fact that $\operatorname*{dom}%
f\subset\operatorname*{dom}(\sup_{k\geq1,\text{ }\lambda_{k}>0}f_{k}^{+}),$ we
have\ $\operatorname*{dom}f^{\ast\ast}\subset\operatorname*{dom}(\sup
_{k\geq1,\text{ }\lambda_{k}>0}f_{k}^{+})^{\ast\ast}$, for all $z\in
\operatorname*{dom}f^{\ast\ast},$ and\
\begin{align*}
\left(
{\textstyle\sum_{k\geq1,\text{ }\lambda_{k}>0}}
\lambda_{k}f_{k}^{+}\right)  ^{\ast\ast}(z) &  =%
{\textstyle\sum_{k\geq1,\text{ }\lambda_{k}>0}}
\lambda_{k}(f_{k}^{+})^{\ast\ast}(z)\\
&  =%
{\textstyle\sum_{k\geq1,\text{ }\lambda_{k}>0}}
\lambda_{k}(f_{k}^{\ast\ast})^{+}(z)\leq%
{\textstyle\sum_{k\geq1}}
\lambda_{k}(f_{k}^{\ast\ast})^{+}(z),
\end{align*}
where we used the fact that $0\leq(f_{k}^{+})^{\ast\ast}=(f_{k}^{\ast\ast
})^{+}$ by\ (\ref{maxbicon}). Combining (\ref{re0}), (\ref{re}) and the last
relation, and taking into account that $\overline{\operatorname*{dom}%
f}^{w^{\ast\ast}}=\overline{\operatorname*{dom}f^{\ast\ast}}^{w^{\ast\ast}},$
we conclude that
\begin{align*}
v(\mathcal{P}) &  \leq\inf_{X^{\ast\ast}\text{ }}(f_{0}^{\ast\ast}+%
{\textstyle\sum_{k\geq1}}
\lambda_{k}(f_{k}^{\ast\ast})^{+}+\lambda_{\infty}f_{\infty}^{\ast\ast
}+\mathrm{I}_{\overline{\operatorname*{dom}f}^{w^{\ast\ast}}})\\
& \\
&  \leq\inf_{_{\substack{f_{k}^{\ast\ast}(z)\leq0,\text{ }k\geq1,\text{
}\\f_{\infty}^{\ast\ast}(z)\leq0,\text{ \ }z\in\overline{\operatorname*{dom}%
f}^{w^{\ast\ast}}}}}(f_{0}^{\ast\ast}(z)+%
{\textstyle\sum_{k\geq1}}
\lambda_{k}(f_{k}^{\ast\ast})^{+}(z)+\lambda_{\infty}f_{\infty}^{\ast\ast
}(z)+\mathrm{I}_{\overline{\operatorname*{dom}f}^{w^{\ast\ast}}}(z)),
\end{align*}
which yields\
\[
v(\mathcal{P})\leq\inf_{_{\substack{f_{k}^{\ast\ast}(z)\leq0,\text{ }%
k\geq1,\\f_{\infty}^{\ast\ast}(z)\leq0,\text{ }z\in\overline
{\operatorname*{dom}f^{\ast\ast}}^{w^{\ast\ast}}}}}f_{0}^{\ast\ast
}(z)=v(\mathcal{P}_{\infty}^{\ast\ast}).
\]

\end{dem}

Other related biconjugate-type relaxations\ can also be considered, providing
alternatives to the use of the function $f_{\infty}$ and the set
$\overline{\operatorname*{dom}f^{\ast\ast}}^{w^{\ast\ast}}.$

\begin{cor}
\label{coro}Under the assumptions of Theorem \ref{main1}, we have
\[
v(\mathcal{P}_{1}^{\ast\ast})=v(\mathcal{P}_{2}^{\ast\ast})=v(\mathcal{P}),
\]
where
\[
(\mathcal{P}_{1}^{\ast\ast})\left\{
\begin{array}
[c]{l}%
\inf~~f_{0}^{\ast\ast}(z)\\
\text{s.t. }f_{k}^{\ast\ast}(z)\leq0,\text{ }k\geq1,\text{ }\\
\text{ \ \ \ \ }f_{\infty}^{\ast\ast}(z)\leq0,\\
\text{ \ \ \ \ }z\in\operatorname*{dom}f^{\ast\ast},
\end{array}
\right.  \text{ and }(\mathcal{P}_{2}^{\ast\ast})\left\{
\begin{array}
[c]{l}%
\inf~~f_{0}^{\ast\ast}(z)\\
\text{s.t. }f_{k}^{\ast\ast}(z)\leq0,\text{ }k\geq1,\text{ }\\
\text{ \ \ \ \ }z\in\overline{\operatorname*{dom}(%
{\textstyle\sum_{k\geq1}}
f_{k}^{+})}^{w^{\ast\ast}}.
\end{array}
\right.
\]

\end{cor}

\begin{dem}
First we\textbf{ }establish\ the inequalities
\begin{equation}
v(\mathcal{P}_{\infty}^{\ast\ast})\leq v(\mathcal{P}_{1}^{\ast\ast})\leq
v(\mathcal{P})\text{ and }v(\mathcal{P}_{\infty}^{\ast\ast})\leq
v(\mathcal{P}_{2}^{\ast\ast})\leq v(\mathcal{P}), \label{ser}%
\end{equation}
where $(\mathcal{P}_{\infty}^{\ast\ast})$ is the relaxation in (\ref{pdual});
afterwards, Theorem \ref{main1} will provide the desired equalities.

The inequality $v(\mathcal{P}_{\infty}^{\ast\ast})\leq v(\mathcal{P}_{1}%
^{\ast\ast})$ is straightforward. To show the inequality $v(\mathcal{P}%
_{\infty}^{\ast\ast})\leq v(\mathcal{P}_{2}^{\ast\ast})$, we first observe
that every $x\in\operatorname*{dom}(%
{\textstyle\sum_{k\geq1}}
f_{k}^{+})$ satisfies $\lim_{k\rightarrow\infty}f_{k}^{+}(x)=0.$ Hence,
$x\in\operatorname*{dom}f$ and $f_{\infty}(x)\leq\limsup_{k\rightarrow\infty
}f_{k}^{+}(x)=0.$ Consequently, $\operatorname*{dom}(%
{\textstyle\sum_{k\geq1}}
f_{k}^{+})\subset\operatorname*{dom}f$ and $\operatorname*{dom}(%
{\textstyle\sum_{k\geq1}}
f_{k}^{+})\subset\lbrack f_{\infty}\leq0].$ Taking $w^{\ast\ast}$-closures
yields
\[
\overline{\operatorname*{dom}(%
{\textstyle\sum_{k\geq1}}
f_{k}^{+})}^{w^{\ast\ast}}\subset\overline{\operatorname*{dom}f}^{w^{\ast\ast
}}=\overline{\operatorname*{dom}f^{\ast\ast}}^{w^{\ast\ast}}\text{, }%
\]
and
\[
\overline{\operatorname*{dom}(%
{\textstyle\sum_{k\geq1}}
f_{k}^{+})}^{w^{\ast\ast}}\subset\overline{[f_{\infty}\leq0]}^{w^{\ast\ast}%
}\subset\overline{[f_{\infty}^{\ast\ast}\leq0]}^{w^{\ast\ast}}=[f_{\infty
}^{\ast\ast}\leq0].
\]
In other words, $F(\mathcal{P}_{2}^{\ast\ast})\subset F(\mathcal{P}_{\infty
}^{\ast\ast}),$ and it follows that $v(\mathcal{P}_{\infty}^{\ast\ast})\leq
v(\mathcal{P}_{2}^{\ast\ast}).$

To prove the inequality $v(\mathcal{P}_{1}^{\ast\ast})\leq v(\mathcal{P}),$ we
take a feasible point $x\in F(\mathcal{P}).$ Then
\[
f_{k}^{\ast\ast}(x)\leq f_{k}(x)\leq0,\text{ for all }k\geq1,
\]
so\ that $x\in\operatorname*{dom}f\subset\operatorname*{dom}f^{\ast\ast}$ and
$f_{\infty}(x)\leq0.$ Thus,
\[
v(\mathcal{P}_{1}^{\ast\ast})\leq\inf_{\substack{f_{k}^{\ast\ast}%
(x)\leq0,\text{ }k\geq1,\text{ }f_{\infty}^{\ast\ast}(x)\leq0\text{ }%
\\x\in\operatorname*{dom}f^{\ast\ast},\text{ }x\in X}}f_{0}^{\ast\ast}%
(x)\leq\inf_{f_{k}(x)\leq0,\text{ }k\geq1\text{ }}f_{0}(x);
\]
that is, $v(\mathcal{P}_{1}^{\ast\ast})\leq v(\mathcal{P}).$ Similarly, every
feasible point $x\in F(\mathcal{P}),$ satisfies\ $%
{\textstyle\sum_{k\geq1}}
f_{k}^{+}(x)\leq0.$ So, $F(\mathcal{P})\subset\operatorname*{dom}(%
{\textstyle\sum_{k\geq1}}
f_{k}^{+})$ and we conclude that
\[
F(\mathcal{P})\subset\overline{F(\mathcal{P})}^{w^{\ast\ast}}\subset
\overline{\operatorname*{dom}(%
{\textstyle\sum_{k\geq1}}
f_{k}^{+})}^{w^{\ast\ast}}.
\]
This shows that $v(\mathcal{P}_{2}^{\ast\ast})\leq v(\mathcal{P}),$ and we are
done.\smallskip
\end{dem}

We consider now a third relaxation:
\[
(\mathcal{P}_{3}^{\ast\ast})\text{ }\left\{
\begin{array}
[c]{l}%
\inf~~f_{0}^{\ast\ast}(z)\\
\text{s.t. }f_{k}^{\ast\ast}(z)\leq0,\text{ }k\geq1,\text{ }\\
\text{ \ \ \ \ }z\in\overline{\operatorname*{dom}f^{\ast\ast}}^{w^{\ast\ast}}.
\end{array}
\right.
\]
The following corollary involves a reinforced Slater condition requiring the
existence of some $x_{0}\in\operatorname*{dom}f_{0}$ and $(\alpha_{k})_{k}%
\in\ell_{1}$ such that $\alpha_{k}>0,$ for all $k\geq1,$ and
\begin{equation}
\sup_{k \ge 1}\alpha_{k}f_{k}(x_{0})<0.\label{ins}%
\end{equation}
Observe that when problem $(\mathcal{P})$ has a finite number of constraints,
this property is equivalent to\ the usual Slater condition. In linear
semi-infinite optimization, (\ref{ins}) is implied by the strong Slater
condition for the so-called reinforced system of constraints introduced in
\cite{GLT97} and \cite[Section 3]{CLP02}.

\begin{cor}
\label{hamdd} Assume that the reinforced Slater condition holds. If
$f:=\sup_{k\geq1}f_{k}$ is continuous somewhere in $\operatorname*{dom}f_{0}$,
then\
\[
v(\mathcal{P}_{3}^{\ast\ast})=v(\mathcal{P}).
\]

\end{cor}

\begin{dem}
Let\ $x_{0}\in\operatorname*{dom}f_{0}$ and $(\alpha_{k})_{k}\in\ell_{1}^{+}$
be as (\ref{ins}). Observe that\ $(\mathcal{P})$ has the same optimal value as
the problem
\[
(\mathcal{P}^{\prime})\left\{
\begin{array}
[c]{l}%
\inf\text{ }f_{0}(x)\\
\text{s.t. }g_{k}(x):=\alpha_{k}f_{k}(x)\leq0,\text{ }k\geq1.
\end{array}
\right.
\]
Since the latter problem satisfies the (usual) Slater condition by assumption,
and the function $\sup_{k\geq1}g_{k}$ is continuous somewhere in
$\operatorname*{dom}f_{0}$, Corollary \ref{coro} entails
\begin{equation}
v(\mathcal{P})=v(\mathcal{P}^{^{\prime}})=\inf_{g_{k}^{\ast\ast}%
(z)\leq0,\text{ }k\geq1,\text{ }z\in\overline{\operatorname*{dom}\left(
{\textstyle\sum_{k\geq1}}
g_{k}^{+}\right)  }^{w^{\ast\ast}}}~~f_{0}^{\ast\ast}(z).\label{tc}%
\end{equation}
Observe that $g_{k}^{+}=\alpha_{k}f_{k}^{+}$ and $g_{k}^{\ast\ast}=\alpha
_{k}f_{k}^{\ast\ast},$ for all $k\geq1.$ Let us show that
\begin{equation}
\overline{\operatorname*{dom}f^{\ast\ast}}^{w^{\ast\ast}}\subset
\overline{\operatorname*{dom}\left(
{\textstyle\sum_{k\geq1}}
g_{k}^{+}\right)  }^{w^{\ast\ast}}.\label{sh}%
\end{equation}
Indeed, given any $x\in\operatorname*{dom}f,$ we obtain
\[%
{\textstyle\sum_{k\geq1}}
g_{k}^{+}(x)\leq%
{\textstyle\sum_{k\geq1}}
\alpha_{k}f_{k}^{+}(x)\leq\left(
{\textstyle\sum_{k\geq1}}
\alpha_{k}\right)  f^{+}(x)<+\infty,
\]
showing that $\operatorname*{dom}f\subset\operatorname*{dom}\left(
{\textstyle\sum_{k\geq1}}
g_{k}^{+}\right)  .$ Hence, taking the $w^{\ast\ast}$-closure yields
\[
\overline{\operatorname*{dom}f^{\ast\ast}}^{w^{\ast\ast}}=\overline
{\operatorname*{dom}f}^{w^{\ast\ast}}\subset\overline{\operatorname*{dom}%
\left(
{\textstyle\sum_{k\geq1}}
g_{k}^{+}\right)  }^{w^{\ast\ast}},
\]
and\ (\ref{tc}) implies
\[
v(\mathcal{P})\leq\inf_{g_{k}^{\ast\ast}(z)\leq0,\text{ }k\geq1,\text{ }%
z\in\overline{\operatorname*{dom}f^{\ast\ast}}^{w^{\ast\ast}}\text{ }}%
~~f_{0}^{\ast\ast}(z)=v(\mathcal{P}_{3}^{\ast\ast}).
\]
Thus, we are done because the inequality\ $v(\mathcal{P}_{3}^{\ast\ast})\leq
v(\mathcal{P})$ is straightforward.
\end{dem}

As consequence of the previous results we obtain the formula for the
biconjugate of the supremum function, which is used later on in section
\ref{sec4}.

\begin{theo}
\label{corsup} Given functions $f_{k}\in\Gamma_{0}(X),$ $k\geq1,$\ we assume
that $f:=\sup_{k\geq1}f_{k}$ is continuous somewhere. Then
\begin{equation}
f^{\ast\ast}=\max\left\{  \sup_{k\geq1}f_{k}^{\ast\ast},\text{ }f_{\infty
}^{\ast\ast}\right\}  +\mathrm{I}_{\overline{\operatorname*{dom}f}%
^{w^{\ast\ast}}}. \label{i1}%
\end{equation}

\end{theo}

\begin{dem}
Assume first that $\inf_{x\in X}f(x)<+\infty.$ Let us write
\begin{equation}
\inf_{x\in X}f(x)=\inf_{f_{k}(x)-y\leq0,\text{ }k\geq1,\text{ }x\in X,\text{
}y\in\mathbb{R}\text{ }}y\equiv\inf_{g_{k}(x,y)\leq0,\text{ }k\geq1,\text{
}x\in X,\text{ }y\in\mathbb{R}\text{ }}g_{0}(x,y),\label{pcor}%
\end{equation}
where $g_{k}\in\Gamma_{0}(X\times\mathbb{R}),$ $k\geq0,$ are defined by
\[
g_{0}(x,y):=y\text{ and }g_{k}(x,y):=f_{k}(x)-y,\text{ for }k\geq1.
\]
We denote $g:=\sup_{k\geq1}g_{k}.$ Then $g(x,y)=\sup_{k\geq1}(f_{k}%
(x)-y)=f(x)-y,$ for all $(x,y)\in X\times\mathbb{R}.$ For all $(z,y)\in
X^{\ast\ast}\times\mathbb{R},$ by (\ref{moreau}) we have%
\[
g^{\ast\ast}(z,y)=\liminf_{x\rightharpoonup z,\text{ }u\rightarrow y,\text{
}x\in X,\text{ }y\in\mathbb{R}\text{ }}g(x,u),
\]
and%
\[
g_{k}^{\ast\ast}(z,y)=\liminf_{x\rightharpoonup z,\text{ }u\rightarrow
y,\text{ }x\in X,\text{ }y\in\mathbb{R}\text{ }}g_{k}(x,u),\text{ for all
}k\geq1;
\]
remember that \textquotedblleft$\rightharpoonup$\textquotedblright\ stands for
the convergence in\ the weak-* topology of\ $X^{\ast\ast}$. Hence,
$g^{\ast\ast}(z,y)=f^{\ast\ast}(z)-y$,\
\[
g_{0}^{\ast\ast}(z,y)=y,\text{ and }g_{k}^{\ast\ast}(z,y)=f_{k}^{\ast\ast
}(z)-y,\text{ for all }k\geq1.
\]
Observe that $g$ is continuous at every point of the form $(x,y)$ with $x$
being a continuity point of $f.$ Moreover, any point $x_{0}\in
\operatorname*{dom}f$ and any $y_{0}>f(x_{0})$ satisfies
\[
g(x_{0},y_{0})=f(x_{0})-y_{0}<0,
\]
showing that $(x_{0},y_{0})$ is a Slater point of the problem in (\ref{pcor})
with the functions $g_{k},$ $k\geq0.$ Therefore, Theorem \ref{main1} yields
\begin{align*}
\inf_{g_{k}(x,y)\leq0,\text{ }k\geq1,\text{ }x\in X,\text{ }y\in
\mathbb{R}\text{ }}g_{0}(x,y) &  =\inf_{\substack{f_{k}^{\ast\ast}%
(z)-y\leq0,\text{ }k\geq1,\text{ }z\in X^{\ast\ast},\text{ }y\in
\mathbb{R}\text{ }\\f_{\infty}^{\ast\ast}(z)-y\leq0,\text{ }z\in
\overline{\operatorname*{dom}(f^{\ast\ast})}^{w^{\ast\ast}}}}y\\
&  =\inf_{z\in X^{\ast\ast}\text{ }}\sup_{k\geq1\text{ }}\{f_{k}^{\ast\ast
}(z),f_{\infty}^{\ast\ast}(z)\}+\mathrm{I}_{\overline{\operatorname*{dom}%
(f^{\ast\ast})}^{w^{\ast\ast}}}(z),
\end{align*}
showing that
\begin{equation}
\inf_{x\in X\text{ }}f(x)=\inf_{z\in X^{\ast\ast}\text{ }}\sup_{k\geq1\text{
}}\{f_{k}^{\ast\ast}(z),f_{\infty}^{\ast\ast}(z)\}+\mathrm{I}_{\overline
{\operatorname*{dom}(f^{\ast\ast})}^{w^{\ast\ast}}}(z).\label{hamd}%
\end{equation}
We now fix $x^{\ast}\in X^{\ast}.$ Since $f\in\Gamma_{0}(X)$, $f^{\ast}$ is
proper and we have
\[
-f^{\ast}(x^{\ast})=\inf_{x\in X\text{ }}\sup_{k\geq1\text{ }}(f_{k}%
(x)-\left\langle x^{\ast},x\right\rangle )<+\infty.
\]
Set $\tilde{f}_{k}:=f_{k}(\cdot)-\left\langle x^{\ast},\cdot\right\rangle ,$
$k\geq1.$ Then $\tilde{f}_{k}\in\Gamma_{0}(X),$ for all $k\geq1,$ and
$\tilde{f}:=\sup_{k\geq1}\tilde{f}_{k}=f(\cdot)-\left\langle x^{\ast}%
,\cdot\right\rangle $ is continuous somewhere. Moreover, using (\ref{moreau}),
we have\ $\tilde{f}_{k}^{\ast\ast}=f_{k}^{\ast\ast}(\cdot)-\left\langle
x^{\ast},\cdot\right\rangle $ and $\limsup_{k}\tilde{f}_{k}=f_{\infty}%
(\cdot)-\left\langle x^{\ast},\cdot\right\rangle $. Next, applying
(\ref{hamd}) with $f_{k}(\cdot)-\left\langle x^{\ast},\cdot\right\rangle $
instead of $f_{k}$ we get\
\[
-f^{\ast}(x^{\ast})=\inf_{z\in X^{\ast\ast}\text{ }}\left(  \sup_{k\geq1\text{
}}\{f_{k}^{\ast\ast}(z),f_{\infty}^{\ast\ast}(z)\}+\mathrm{I}_{\overline
{\operatorname*{dom}(f^{\ast\ast})}^{w^{\ast\ast}}}(z)-\left\langle x^{\ast
},z\right\rangle \right)  ,
\]
which leads us to
\[
f^{\ast}(x^{\ast})=\left(  \sup_{k\geq1\text{ }}\{f_{k}^{\ast\ast},f_{\infty
}^{\ast\ast}\}+\mathrm{I}_{\overline{\operatorname*{dom}(f^{\ast\ast}%
)}^{w^{\ast\ast}}}\right)  ^{\ast}(x^{\ast}),\text{ for all }x^{\ast}\in
X^{\ast}.
\]
Thus, taking the conjugate on both sides and using again (\ref{moreau}), we
obtain
\[
f^{\ast\ast}=\left(  \sup_{k\geq1\text{ }}\{f_{k}^{\ast\ast},f_{\infty}%
^{\ast\ast}\}+\mathrm{I}_{\overline{\operatorname*{dom}(f^{\ast\ast}%
)}^{w^{\ast\ast}}}\right)  ^{\ast\ast}=\sup_{k\geq1\text{ }}\{f_{k}^{\ast\ast
},f_{\infty}^{\ast\ast}\}+\mathrm{I}_{\overline{\operatorname*{dom}%
(f^{\ast\ast})}^{w^{\ast\ast}}},
\]
and the proof is complete.
\end{dem}

\section{Concave-like setting\label{sec3}}

This section focuses on the case where the family of constraint functions is concave-like.
\begin{definiti}
A family $\{f_{t},$ $t\in T\}\subset\Gamma_{0}(X)$ is said to be
\emph{concave-like} if, for every $t_{1},\cdots,t_{m}\in T$ and $\lambda
_{1},\cdots,\lambda_{m}\in\mathbb{R}_{+},$ $m\geq1,$ there exists $t\in T$
such that
\[%
{\textstyle\sum_{k=1}^{m}}
\lambda_{k}f_{t_{k}}\leq(%
{\textstyle\sum_{k=1}^{m}}
\lambda_{k})f_{t}.
\]

\end{definiti}

Due to the convention $0\cdot(+\infty)=+\infty,$ the concave-like property is
equivalent to the closedness for convex combinations (see \cite[Definition
5.1.3]{CHLBook}).

The following remark highlights two important features of concave-like
settings. 

\begin{rem}
Let $\{f_{t},$ $t\in T\}\subset\Gamma_{0}(X)$ be a concave-like family and
denote $f:=\sup_{t\in T}f_{t}$. Assume that\ $T$ is compact and that the
mappings $t\mapsto f_{t}(x),$ $x\in\operatorname*{dom}f,$ are upper semicontinuous.

$(i)$ If $v(\mathcal{P})$ is finite and $(\mathcal{P})$ satisfies the Slater
condition, then Lemma \ref{lemHKL} takes the following specific form,
\[
v(\mathcal{P})=\max_{t\in T,\text{ }\lambda\geq0\text{ }}\inf_{x\in
\operatorname*{dom}f\text{ }}(f_{0}(x)+\lambda f_{t}(x)).
\]
A proof of this result can be obtained by following the arguments in
\cite[Corollary 6]{CHL26-JOTA}.

$(ii)$ According to \cite[Corollary 6]{CHL26-JOTA}, the following biconjugate
supremum rule holds,%
\begin{equation}
(\sup_{t\in T}f_{t})^{\ast\ast}=\sup_{t\in T}(f_{t}+\mathrm{I}%
_{\operatorname*{dom}f})^{\ast\ast}.\label{ghaf}%
\end{equation}

\end{rem}

In this concave-like setting, we introduce the following relaxed problem:
\begin{equation}
(\mathcal{P}_{c}^{\ast\ast})\text{ }\left\{
\begin{array}
[c]{l}%
\inf~~(f_{0}+\mathrm{I}_{\operatorname*{dom}f})^{\ast\ast}(z)\\
\text{s.t. }(f_{t}+\mathrm{I}_{(\operatorname*{dom}f_{0})\cap
(\operatorname*{dom}f)})^{\ast\ast}(z)\leq0,\text{ }t\in T,\text{ }z\in
X^{\ast\ast}.
\end{array}
\right.  \label{tal1}%
\end{equation}
Moreover, when $f$ is continuous somewhere in $\operatorname*{dom}f_{0}$, we
also prove that the last problem reduces to
\begin{equation}
\left\{
\begin{array}
[c]{l}%
\inf~~f_{0}^{\ast\ast}(z)\\
\text{s.t. }f_{t}^{\ast\ast}(z)\leq0,\text{ }t\in T,\text{ }\\
\text{ \ \ \ \ }z\in\overline{\operatorname*{dom}f}^{w^{\ast\ast}}.
\end{array}
\right.  \label{tal2}%
\end{equation}

\begin{theo}
\label{marco11}Given a concave-like family $\{f_{t},$ $t\in T\}\subset
\Gamma_{0}(X)$, we assume that $T$ is compact and the mappings $t\mapsto
f_{t}(x),$ $x\in\operatorname*{dom}f$ (assumed nonempty), are upper
semicontinuous. Then, provided that $(\mathcal{P})$ satisfies the Slater
condition, we have
\begin{equation}
v(\mathcal{P}_{c}^{\ast\ast})=v(\mathcal{P}). \label{b1}%
\end{equation}
When, additionally, $f$ is continuous somewhere in $\operatorname*{dom}f_{0}$,
(\ref{b1}) remains valid if $(\mathcal{P}_{c}^{\ast\ast})$ stands for the
problem in (\ref{tal2}).
\end{theo}

\begin{dem}
Notice first that
\begin{align*}
v(\mathcal{P}_{c}^{\ast\ast})  &  \leq\inf_{(f_{t}+\mathrm{I}%
_{(\operatorname*{dom}f_{0})\cap(\operatorname*{dom}f)})^{\ast\ast}%
(x)\leq0,\text{ }x\in X,\text{ }t\in T\text{ }}\text{ }(f_{0}+\mathrm{I}%
_{\operatorname*{dom}f})^{\ast\ast}(x)\\
&  \leq\inf_{f_{t}(x)\leq0,\text{ }x\in(\operatorname*{dom}f_{0}%
)\cap(\operatorname*{dom}f),\text{ }t\in T\text{ }}\text{ }f_{0}%
(x)+\mathrm{I}_{\operatorname*{dom}f}(x)=v(\mathcal{P}).
\end{align*}
So, $v(\mathcal{P}_{c}^{\ast\ast})\leq v(\mathcal{P})$ and we may assume that
$v(\mathcal{P})\in\mathbb{R};$ in fact, if $v(\mathcal{P})=-\infty,$ then the
weak inequality implies $v(\mathcal{P}_{c}^{\ast\ast})=v(\mathcal{P})=-\infty
$, and we are done.

To proceed, we write\
\begin{align*}
0 &  =\inf_{x\in X\text{ }}\sup\{f_{0}(x)-v(\mathcal{P});\text{ }%
f_{t}(x),\text{ }t\in T\}\\
&  =\inf_{z\in X^{\ast\ast}\text{ }}(\sup\{f_{0}-v(\mathcal{P});\text{ }%
f_{t},\text{ }t\in T\})^{\ast\ast}(z),
\end{align*}
because a function and its biconjugate have the same infimum. Thus, applying
(\ref{b1}) to the compact set $T\cup\{0\}$, we obtain
\begin{equation}
0=\inf_{z\in X^{\ast\ast}\text{ }}\sup\{(f_{0}+\mathrm{I}_{\operatorname*{dom}%
f})^{\ast\ast}(z)-v(\mathcal{P});\text{ }(f_{t}+\mathrm{I}%
_{(\operatorname*{dom}f_{0})\cap(\operatorname*{dom}f)})^{\ast\ast}(z),\text{
}t\in T\}.\label{oi}%
\end{equation}
Moreover, since $(\mathcal{P})$ satisfies the Slater condition, there exists
$x_{0}\in\operatorname*{dom}f_{0}$ such that $f(x_{0})<0.$ In particular,
$x_{0}\in(\operatorname*{dom}f_{0})\cap(\operatorname*{dom}f)$ and
\[
(f_{0}+\mathrm{I}_{\operatorname*{dom}f})^{\ast\ast}(x_{0})\leq(f_{0}%
+\mathrm{I}_{\operatorname*{dom}f})(x_{0})=f_{0}(x_{0})<+\infty;
\]
that is, $x_{0}\in\operatorname*{dom}(f_{0}+\mathrm{I}_{\operatorname*{dom}%
f})^{\ast\ast}.$ Furthermore, we have
\[
\sup_{t\in T\text{ }}(f_{t}+\mathrm{I}_{(\operatorname*{dom}f_{0}%
)\cap(\operatorname*{dom}f)})^{\ast\ast}(x_{0})\leq\sup_{t\in T\text{ }}%
(f_{t}+\mathrm{I}_{(\operatorname*{dom}f_{0})\cap(\operatorname*{dom}%
f)})(x_{0})\leq f(x_{0})<0,
\]
and therefore $x_{0}$ is also a Slater point of $(\mathcal{P}_{c}^{\ast\ast
}).$ Consequently, by (\ref{oi}), relation (\ref{b1}) follows by using
(\ref{ca}).

Finally, assume that\ $f$ is continuous at some point in $\operatorname*{dom}%
f_{0}$. Then, by the biconjugate sum rule in (\ref{sumbicon}), we have
\[
(f_{0}+\mathrm{I}_{\operatorname*{dom}f})^{\ast\ast}=f_{0}^{\ast\ast
}+\mathrm{I}_{\operatorname*{dom}f}^{\ast\ast}=f_{0}^{\ast\ast}+\mathrm{I}%
_{\overline{\operatorname*{dom}f}^{w^{\ast\ast}}}%
\]
and, for all $t\in T,$
\[
(f_{t}+\mathrm{I}_{(\operatorname*{dom}f_{0})\cap(\operatorname*{dom}%
f)})^{\ast\ast}=f_{t}^{\ast\ast}+\mathrm{I}_{\overline{(\operatorname*{dom}%
f_{0})\cap(\operatorname*{dom}f)}^{w^{\ast\ast}}}=f_{t}^{\ast\ast}%
+\mathrm{I}_{\overline{\operatorname*{dom}f}^{w^{\ast\ast}}}+\mathrm{I}%
_{\overline{\operatorname*{dom}f_{0}}^{w^{\ast\ast}}}.
\]
At the same time, relation (\ref{hamdo}) ensures that\ $\overline
{(\operatorname*{dom}f_{0})\cap(\operatorname*{dom}f)}^{w^{\ast\ast}%
}=\overline{\operatorname*{dom}f_{0}}^{w^{\ast\ast}}\cap\overline
{\operatorname*{dom}f}^{w^{\ast\ast}}.$ Consequently,
\begin{align*}
v(\mathcal{P}_{c}^{\ast\ast}) &  =\inf_{(f_{t}+\mathrm{I}%
_{(\operatorname*{dom}f_{0})\cap(\operatorname*{dom}f)})^{\ast\ast}%
(z)\leq0,\text{ }t\in T\text{ }}\text{ }(f_{0}+\mathrm{I}_{\operatorname*{dom}%
f})^{\ast\ast}(z)\\
&  =\inf_{f_{t}^{\ast\ast}(z)\leq0,\text{ }t\in T\text{ , }z\in\overline
{\operatorname*{dom}f_{0}}^{w^{\ast\ast}}\cap\overline{\operatorname*{dom}%
f}^{w^{\ast\ast}}\text{ }}\text{ }f_{0}^{\ast\ast}(z)+\mathrm{I}%
_{\overline{\operatorname*{dom}f}^{w^{\ast\ast}}}(z)\\
&  =\inf_{f_{t}^{\ast\ast}(z)\leq0,\text{ },\text{ }t\in T\text{, }%
z\in\overline{\operatorname*{dom}f}^{w^{\ast\ast}}}\text{ }f_{0}^{\ast\ast
}(z),
\end{align*}
and thus (\ref{b1}) is satisfied when $(\mathcal{P}_{c}^{\ast\ast})$ is
defined as in (\ref{tal2}).
\end{dem}

\section{Relationship with duality in convex programming\label{sec4a}}

In this section we emphasize the connection between the biconjugate
relaxation, studied in the previous sections, and the Fenchel duality in
infinite convex optimization.

Consider the (primal) convex optimization problem\
\[
(\mathcal{P})\text{ \ \ }\inf_{X\text{ }}(f+g)(x),
\]
where\ $f,$ $g:X\rightarrow\mathbb{R}_{\infty}$ are two proper convex
functions. The associated (Fenchel)\ \emph{dual} problem is given by
\begin{equation}
(\mathcal{D})\text{ \ }\sup_{x^{\ast}\in X^{\ast}\text{ }}(-f^{\ast}(x^{\ast
})-g^{\ast}(-x^{\ast})).\label{Fencheld}%
\end{equation}
We refer to\ $(\mathcal{P},\mathcal{D})$ as a dual pair and assume that
$v(\mathcal{P})<+\infty.$ The classical strong duality theorem (see, e.g.,
\cite{Mo66}) asserts that if either $f$ or $g$ is finite and continuous at
some point in the effective domain of the other function, then\
\begin{equation}
v(\mathcal{P})=v(\mathcal{D}).\label{salto}%
\end{equation}
Moreover, problem $(\mathcal{D})$ admits optimal solutions, leading to the
so-called strong duality between $(\mathcal{P})$ and $(\mathcal{D}%
)$.\ Furthermore, treating\ $(\mathcal{D})$ itself as a primal problem, its
Fenchel dual turns out to be the biconjugate relaxation of $(\mathcal{P});$
namely,
\[
(\mathcal{D}^{\ast})\text{\ \ }\inf_{z\in X^{\ast\ast}\text{ }}(f^{\ast\ast
}(z)+g^{\ast\ast}(z)).
\]
Applying\ the weak duality theorem to the dual pairs $(\mathcal{P}%
,\mathcal{D})$ and $(\mathcal{D},\mathcal{D}^{\ast}),$ we obtain
\begin{equation}
v(\mathcal{D}^{\ast})\leq v(\mathcal{D})\leq v(\mathcal{P}).\label{sandwish}%
\end{equation}
Therefore, having zero-duality gaps for the pairs $(\mathcal{P},\mathcal{D})$
and $(\mathcal{D},\mathcal{D}^{\ast})$ is equivalent to a zero gap between
$(\mathcal{P})$ and $(\mathcal{D}^{\ast}).$

The following preliminary result in this section shows that the optimal values
of $(\mathcal{P})$ and $(\mathcal{D}^{\ast})$\ coincide under the same
continuity condition imposed above only on the original functions $f$ and $g.$
By the Moreau--Rockafellar theorem (\ref{moreau}), this issue is relevant only
in nonreflexive Banach spaces.

\begin{prop}
\label{prop}Let $f$ and $g$ be convex functions. Assume that one of them is
finite and continuous at some point in the effective domain of the other.
Then
\[
v(\mathcal{D}^{\ast})=v(\mathcal{P}).
\]
Consequently, there is a zero-duality gap between $(\mathcal{P})$ and
$(\mathcal{D}).$
\end{prop}

\begin{dem}
Since $v(\mathcal{D}^{\ast})\leq v(\mathcal{P})$, we may assume that
$v(\mathcal{P})\in\mathbb{R}$. According to\ \cite{HLZ08}, we have
$\operatorname*{cl}(f+g)=(\operatorname*{cl}f)+(\operatorname*{cl}g).$ Thus,
since
\[
v(\mathcal{P})=\inf_{X}(f+g)=\inf_{X}\operatorname*{cl}(f+g)=\inf
_{X}((\operatorname*{cl}f)+(\operatorname*{cl}g)),
\]
it follows that $\operatorname*{cl}f,$ $\operatorname*{cl}g\in\Gamma_{0}(X).$
Applying (\ref{sumbicon}) to the functions $(\operatorname*{cl}f)$ and
$(\operatorname*{cl}g)$,\ the current assumption\ ensures that
\[
(f+g)^{\ast\ast}=(\operatorname*{cl}(f+g))^{\ast\ast}=((\operatorname*{cl}%
f)+(\operatorname*{cl}g))^{\ast\ast}=(\operatorname*{cl}f)^{\ast\ast
}+(\operatorname*{cl}g)^{\ast\ast}=f^{\ast\ast}+g^{\ast\ast}.
\]
Thus,
\[
v(\mathcal{P})=\inf_{X}(f+g)=\inf_{X^{\ast\ast}}(f+g)^{\ast\ast}=\inf
_{X^{\ast\ast}}(f^{\ast\ast}+g^{\ast\ast})=v(\mathcal{D}^{\ast}),
\]
and the first conclusion holds. The second assertion follows
from\ (\ref{sandwish}).
\end{dem}

The crucial step in the proof above is the continuity assumption,
which\ ensures that $(f+g)^{\ast\ast}=f^{\ast\ast}+g^{\ast\ast}.$ As the
following example shows, this hypothesis cannot be dropped.

\begin{exam}
Let $X=c_{0},$ and let $L\subset c_{0}$ be a nontrivial finite-dimensional
subspace such that $L\cap c_{00}=\{\theta\},$ where $c_{00}$ denotes the
subspace of real sequences with only finitely many zeros. Define the convex
functions $f:=\mathrm{I}_{c_{00}}$ and $g:=\mathrm{I}_{L}.$ Then, by
(\ref{moreau}), we have $f^{\ast\ast}=\mathrm{I}_{\operatorname*{cl}%
^{w^{\ast\ast}}(c_{00})}=\mathrm{I}_{\ell_{\infty}}$ and $g^{\ast\ast
}:=\mathrm{I}_{L},$\ so that
\[
f+g=(f+g)^{\ast\ast}=\mathrm{I}_{\{\theta\}}\text{ and }f^{\ast\ast}%
+g^{\ast\ast}=\mathrm{I}_{L}.
\]
Hence, the equality $(f+g)^{\ast\ast}=f^{\ast\ast}+g^{\ast\ast}$ fails.
\end{exam}

We come back to the infinite\ convex program
\begin{equation}
(\mathcal{P})\text{ \ }\left\{  \text{\ }%
\begin{array}
[c]{l}%
\inf~~f_{0}(x)\\
\text{s.t. }f_{k}(x)\leq0,\text{ }k\geq1,\\
\text{ \ \ \ \ }x\in X,
\end{array}
\right.  \label{pdef}%
\end{equation}
where $f_{k}:X\rightarrow\mathbb{R}_{\infty},$ $k\geq1,$ are proper and
convex. We assume\ that $v(\mathcal{P})<+\infty$ and that $(\mathcal{P})$
satisfies the Slater condition. Observe that\
\[
v(\mathcal{P})=\inf_{x\in X\text{ }}\ f_{0}(x)+\mathrm{I}_{[f\leq0]}%
(x)=\inf_{x\in X\text{ }}\ f_{0}(x)+\mathrm{I}_{C}(x),
\]
where\ $f:=\sup_{k\geq1}f_{k}$ and
\[
C:=[f\leq0]\cap\lbrack f_{\infty}\leq0]\cap\operatorname*{dom}f.
\]
The set\ $C$ is nonempty as $v(\mathcal{P})<+\infty.$ Since $(\mathrm{I}%
_{C})^{\ast}=\mathrm{\sigma}_{C},$ the Fenchel dual (see (\ref{Fencheld})) of
$(\mathcal{P})$ is written as
\begin{equation}
(\mathcal{D})\text{ \ }\sup_{x^{\ast}\in X^{\ast}\text{ }}\ -f_{0}^{\ast
}(x^{\ast})-\mathrm{\sigma}_{C}(-x^{\ast}). \label{fd}%
\end{equation}
Since $f\geq f_{k}$ for all $k\geq1,$ $f_{\infty}\geq f_{\infty}^{\ast\ast},$
and $\mathrm{I}_{\operatorname*{dom}f}\geq\mathrm{I}_{\overline
{\operatorname*{dom}f}^{w^{\ast\ast}}},$ the system
\[
\{f_{k}^{\ast\ast}\leq0;k\geq1,f_{\infty}^{\ast\ast}\leq0,\mathrm{I}%
_{\overline{\operatorname*{dom}f}^{w^{\ast\ast}}}\leq0\}
\]
inherits the Slater condition from the system $\{f_{k}\leq0,$ $k\geq1\}.$
Thus, taking into account that
\[
C\subset C^{\prime}:=[\sup_{k\geq1}f_{k}^{\ast\ast}\leq0]\cap\lbrack
f_{\infty}^{\ast\ast}\leq0]\cap\overline{\operatorname*{dom}f}^{w^{\ast\ast}%
},
\]
we get
\[
\mathrm{\sigma}_{C}\leq\mathrm{\sigma}_{C^{\prime}}=\min_{\alpha>0}%
(\alpha(\max\{\sup_{k\geq1}f_{k}^{\ast\ast},f_{\infty}^{\ast\ast}%
,\mathrm{I}_{\overline{\operatorname*{dom}f}^{w^{\ast\ast}}}\}))^{\ast},
\]
where the last equality follows\ from \cite[Theorem 3.3.4]{CHLBook}.
Therefore, the new dual problem
\[
(\mathcal{D}^{\prime})\text{ \ }\sup_{x^{\ast}\in X^{\ast}\text{ }}%
\ -f_{0}^{\ast}(x^{\ast})-\min_{\alpha>0}(\alpha(\max\{\sup\limits_{k\geq
1}f_{k}^{\ast\ast},f_{\infty}^{\ast\ast},\mathrm{I}_{\overline
{\operatorname*{dom}f}^{w^{\ast\ast}}}\})))^{\ast}(-x^{\ast}),
\]
has an optimal value that provides a lower bound\ for the Fenchel dual
$(\mathcal{D})$ in (\ref{fd}).$\ $In contrast to $(\mathcal{D}),$ problem
$(\mathcal{D}^{\prime})$ explicitly involves the conjugates of the data
functions $f_{k},$ $k\geq1.$ Using again \cite[Theorem 3.3.1]{CHLBook}, we see
that the Fenchel dual of $(\mathcal{D}^{\prime})$ coincides with\ the
biconjugate relaxation $(\mathcal{P}_{\infty}^{\ast\ast})$ given in
(\ref{pdual}) by
\[
(\mathcal{P}_{\infty}^{\ast\ast})\text{ }\left\{
\begin{array}
[c]{l}%
\inf~~f_{0}^{\ast\ast}(z)\\
\text{s.t. }f_{k}^{\ast\ast}(z)\leq0,\text{ }k\geq1,\text{ }\\
\text{ \ \ \ \ }f_{\infty}^{\ast\ast}(z)\leq0,\\
\text{ \ \ \ \ }z\in\overline{\operatorname*{dom}(f^{\ast\ast})}^{w^{\ast\ast
}}.
\end{array}
\right.
\]
Then, by the weak duality,\ we have
\begin{equation}
v(\mathcal{P})\geq v(\mathcal{D})\geq v(\mathcal{D}^{^{\prime}})\geq
v(\mathcal{P}_{\infty}^{\ast\ast}). \label{serie}%
\end{equation}
The following theorem completes the announced objective of the section,
showing the strong relation between our biconjugate relaxation and the Fenchel
duality in infinite convex optimization.

\begin{prop}
\label{prop2b}Suppose that $(\mathcal{P})$ in (\ref{pdef}) satisfies the Slater
condition. If the function $f:=\sup_{k\geq1}f_{k}$ is continuous at some point
in $\operatorname*{dom}f_{0},$ then\
\[
v(\mathcal{P})=v(\mathcal{D})=v(\mathcal{D}^{^{\prime}})=v(\mathcal{P}%
_{\infty}^{\ast\ast}).
\]

\end{prop}

\begin{dem}
According to Theorem \ref{main1}, we have $v(\mathcal{P})=v(\mathcal{P}%
_{\infty}^{\ast\ast}).$ Then the conclusion follows by (\ref{serie}).
\end{dem}

\section{Biconjugate relaxation $(\mathcal{P}^{\ast\ast})$ is not
enough\label{sec4}}

We discuss in this section an example exhibiting a gap between problem
$(\mathcal{P})$ and its biconjugate\ relaxation $(\mathcal{P}^{\ast\ast}),$
given in (\ref{eti}) by
\[
(\mathcal{P}^{\ast\ast})\text{ \ }\left\{
\begin{array}
[c]{l}%
\inf~~f_{0}^{\ast\ast}(z)\\
\text{s.t. }f_{k}^{\ast\ast}(z)\leq0,\text{ }k\geq1,\\
\text{ \ \ \ \ }z\in X^{\ast\ast},
\end{array}
\right.
\]
thereby justifying the necessity of considering\ alternative relaxations such
as\ $(\mathcal{P}_{\infty}^{\ast\ast}).$

\begin{exam}
Let\ $c_{0}$ be the Banach space of real sequences converging to zero, endowed
with the supremum norm. Its bidual space is\ $c_{0}^{\ast\ast}=\ell_{\infty}.$
We consider the following optimization problem given in $c_{0}\times
\mathbb{R}$ by
\[
(\mathcal{P})\text{ \ }\left\{
\begin{array}
[c]{l}%
\inf~~f_{0}(x,y):=y+%
{\textstyle\sum_{k\geq1}}
2^{-k}\left\vert x_{k}-1\right\vert \\
\text{s.t. }f_{k}(x,y):=\left\vert x_{k}-1\right\vert -y\leq0,\text{ }%
k\geq1,\\
~\text{\ \ \ \ }x:=(x_{k})_{k}\in c_{0},\text{ }y\in\mathbb{R}.
\end{array}
\right.
\]
We associate to $(\mathcal{P})$ the function\
\[
f(x,y):=\sup_{k\geq1\text{ }}f_{k}(x,y)=\sup_{k\geq1\text{ }}(\left\vert
x_{k}-1\right\vert -y),\text{ }x\in c_{0},\text{ }y\in\mathbb{R}.
\]
We proceed by checking the followings facts:

\begin{itemize}
\item Problem $(\mathcal{P})$ satisfies the Slater condition.

Let $0_{\mathbb{N}}$ be the null sequence in $c_{0}.$ Then $(0_{\mathbb{N}%
},2)\in\operatorname*{dom}f_{0}=c_{0}\times\mathbb{R}$ is a Slater point for
$(\mathcal{P}).$

\item The function $f$ is continuous on $c_{0}\times\mathbb{R}$.

Since $f(x,y)\leq\left\Vert x\right\Vert +\left\vert y\right\vert +1,$ for all
$(x,y)\in c_{0}\times\mathbb{R}$, the (convex) function $f$ is continuous on
$c_{0}\times\mathbb{R}$.

\item We have $f_{k}^{\ast\ast}(z,y)=\left\vert z_{k}-1\right\vert -y,$ for
all $(z,y)\in\ell_{\infty}\times\mathbb{R}.$

Indeed, since $f_{k}\in\Gamma_{0}(c_{0}\times\mathbb{R}),$ it follows from
(\ref{moreau}) that
\[
f_{k}^{\ast\ast}(z,y)=\liminf_{x\rightharpoonup z,\text{ }u\rightarrow
y,\text{ }x\in X\text{ }}(\left\vert x_{k}-1\right\vert -u)=\left\vert
z_{k}-1\right\vert -y.
\]

\item We have $f_{0}^{\ast\ast}(z,y)=y+%
{\textstyle\sum_{k\geq1}}
2^{-k}\left\vert z_{k}-1\right\vert ,$ for all $z\in\ell_{\infty}$ and
$y\in\mathbb{R}.$

To see this, we introduce the convex\ functions $g,$ $g_{k}:c_{0}%
\rightarrow\mathbb{R},$ $k\geq1,$ defined by
\[
g_{k}(x):=\left\vert x_{k}-1\right\vert ,\text{ }g(x):=%
{\textstyle\sum_{k\geq1}}
2^{-k}\left\vert x_{k}-1\right\vert .
\]
Each $g_{k}$ is continuous on $c_{0}.$ In addition, since $g\leq\left\Vert
\cdot\right\Vert +1,$ we have\ $g^{\ast\ast}\leq(\left\Vert \cdot\right\Vert
)^{\ast\ast}+1=\left\Vert \cdot\right\Vert +1$ and, so, $\operatorname*{dom}%
g^{\ast\ast}=\ell_{\infty}.$ Moreover, (\ref{moreau}) implies
\[
g_{k}^{\ast\ast}(z)=\liminf_{x\rightharpoonup z,\text{ }x\in X\text{ }%
}\left\vert x_{k}-1\right\vert =\left\vert z_{k}-1\right\vert ,\text{ for all
}z\in\ell_{\infty}.
\]
Consequently,\ Lemma \ref{prop3} gives rise, for all $z\in\ell_{\infty
}=\operatorname*{dom}g^{\ast\ast},$ to
\[
g^{\ast\ast}(z)=%
{\textstyle\sum_{k\geq1}}
2^{-k}g_{k}^{\ast\ast}(z)=%
{\textstyle\sum_{k\geq1}}
2^{-k}\left\vert z_{k}-1\right\vert .
\]
Therefore, once again by (\ref{moreau}),
\begin{align*}
f_{0}^{\ast\ast}(z,y)  &  =\liminf_{x\rightharpoonup z,\text{ }u\rightarrow
y,\text{ }x\in X\text{ }}(u+%
{\textstyle\sum_{k\geq1}}
2^{-k}\left\vert x_{k}-1\right\vert )\\
&  =y+g^{\ast\ast}(z)=y+%
{\textstyle\sum_{k\geq1}}
2^{-k}\left\vert z_{k}-1\right\vert .
\end{align*}
We are now in position to write the biconjugate relaxation $(\mathcal{P}%
^{\ast\ast}):$%
\[
(\mathcal{P}^{\ast\ast})\text{ }\left\{
\begin{array}
[c]{l}%
\inf~y+%
{\textstyle\sum_{k\geq1}}
2^{-k}\left\vert z_{k}-1\right\vert \\
\text{s.t. }\left\vert z_{k}-1\right\vert -y\leq0,\text{ }k\geq1,\\
\text{ \ \ \ \ }z:=(z_{k})_{k}\in\ell_{\infty},\text{ }y\in\mathbb{R}.
\end{array}
\right.
\]

\item There is a gap between $(\mathcal{P})$ and $(\mathcal{P}^{\ast\ast}).$

Indeed, on the one hand, we have
\[
v(\mathcal{P})=\inf_{x\in c_{0}\text{ }}\left(  \sup_{k\geq1\text{ }%
}\left\vert x_{k}-1\right\vert +%
{\textstyle\sum_{k\geq1}}
2^{-k}\left\vert x_{k}-1\right\vert \right)  .
\]
Consider the sequence $x_{n}=(1,\overset{n}{\cdots},1,0,0,\cdots),$ consisting
of $n$ ones followed by zeros, and define
\[
h(x):=\sup_{k\geq1}\left\vert x_{k}-1\right\vert +%
{\textstyle\sum_{k\geq1}}
2^{-k}\left\vert x_{k}-1\right\vert ,\text{ }x\in c_{0}.
\]
Then,
\[
v(\mathcal{P})\leq\liminf_{n}h(x_{n})=\liminf_{n}(1+%
{\textstyle\sum_{k\geq n+1}}
\frac{1}{2^{k}})=1.
\]
At the same time, $v(\mathcal{P})\geq\inf_{x\in c_{0}}\left(  \sup_{k\geq
1}\left\vert x_{k}-1\right\vert \right)  =1,$ and we get $v(\mathcal{P})=1.$

On the other hand, since $1_{\mathbb{N}}:=(1,1,\cdots)\in\ell_{\infty},$ we
obtain
\[
v(\mathcal{P}^{\ast\ast})=\inf_{x\in\ell_{\infty}\text{ }}\left(  \sup
_{k\geq1\text{ }}\left\vert x_{k}-1\right\vert +%
{\textstyle\sum_{k\geq1}}
2^{-k}\left\vert x_{k}-1\right\vert \right)  \leq0.
\]
But we have $v(\mathcal{P}^{\ast\ast})\geq0$, and so $v(\mathcal{P}^{\ast\ast
})=0<1=v(\mathcal{P}).$ This proves the existence of a duality gap between
$(\mathcal{P})$ and $(\mathcal{P}^{\ast\ast})$.
\end{itemize}

Therefore, the need arises for an alternative biconjugate relaxation
exhibiting a zero duality gap with $(\mathcal{P})$. To this end, we consider
the relaxation $(\mathcal{P}_{\infty}^{\ast\ast})$ given in (\ref{pdual}),
namely
\[
(\mathcal{P}_{\infty}^{\ast\ast})\text{ }\left\{
\begin{array}
[c]{l}%
\inf~~f_{0}^{\ast\ast}(z)\\
\text{s.t. }f_{k}^{\ast\ast}(z)\leq0,\text{ }k\geq1,\text{ }\\
\text{ \ \ \ \ }f_{\infty}^{\ast\ast}(z)\leq0,\\
\text{ \ \ \ \ }z\in\overline{\operatorname*{dom}(f^{\ast\ast})}^{w^{\ast\ast
}}.
\end{array}
\right.
\]
We next compute the function $f_{\infty}^{\ast\ast}$ and identify the set
$\overline{\operatorname*{dom}(f^{\ast\ast})}^{w^{\ast\ast}}.$

\begin{itemize}
\item We have $f_{\infty}(x,y)=1-y,$ for all $(x,y)\in c_{0}\times\mathbb{R}.$

In fact, $f_{\infty}(x,y)=\limsup_{k\rightarrow\infty\text{ }}f_{k}%
(x,y)=\lim_{k\rightarrow\infty}(\left\vert x_{k}-1\right\vert -y)=1-y.$

\item We have $f_{\infty}^{\ast\ast}(z,y)=1-y,$ for all $(z,y)\in\ell_{\infty
}\times\mathbb{R}.$

Since $f_{\infty}\in\Gamma_{0}(c_{0}\times\mathbb{R}),$ (\ref{moreau}) entails
$f_{\infty}^{\ast\ast}(z,y)=\liminf_{x\rightharpoonup z,\text{ }u\rightarrow
y,\text{ }x\in X\text{ }}f_{\infty}(z,u)=\lim_{u\rightarrow y}(1-u)=1-y.$

\item We have $\overline{\operatorname*{dom}f}^{w^{\ast\ast}}=\ell_{\infty
}\times\mathbb{R}.$

We have shown above that\ $\operatorname*{dom}f=c_{0}\times\mathbb{R}.$ So, by
Goldstein's theorem, we obtain $\overline{\operatorname*{dom}f}^{w^{\ast\ast}%
}=\overline{c_{0}\times\mathbb{R}}^{w^{\ast\ast}}=\ell_{\infty}\times
\mathbb{R}.$

Consequently, $(\mathcal{P}_{\infty}^{\ast\ast})$ takes the form
\[
(\mathcal{P}_{\infty}^{\ast\ast})\text{ }\left\{
\begin{array}
[c]{l}%
\inf~~y+%
{\textstyle\sum_{k\geq1}}
2^{-k}\left\vert z_{k}-1\right\vert \\
\text{s.t. }\left\vert z_{k}-1\right\vert -y\leq0,\text{ }k\geq1,\text{ }\\
\text{ \ \ \ \ }1-y\leq0,\\
\text{ \ \ \ \ }z\in\ell_{\infty},\text{ }y\in\mathbb{R}.
\end{array}
\right.
\]

\item There is no gap between $(\mathcal{P})$ and $(\mathcal{P}_{\infty}%
^{\ast\ast})$.

Indeed, this is a consequence of Theorem \ref{main1}. For completeness, we
verify this claim. We have
\[
v(\mathcal{P}_{\infty}^{\ast\ast})=\inf_{z\in\ell_{\infty}}(\max\{\sup
_{k\geq1}\left\vert z_{k}-1\right\vert ,1\}+%
{\textstyle\sum_{k\geq1}}
2^{-k}\left\vert z_{k}-1\right\vert ).
\]
In particular, evaluating the last objective function at $1_{\mathbb{N}}%
\in\ell_{\infty}$ yields
\[
v(\mathcal{P}_{\infty}^{\ast\ast})\leq1.
\]
Hence, since $v(\mathcal{P}_{\infty}^{\ast\ast})\geq1,$ we get $v(\mathcal{P}%
_{\infty}^{\ast\ast})=v(\mathcal{P})=1.$ In other words, there is no gap
between $(\mathcal{P})$ and $(\mathcal{P}_{\infty}^{\ast\ast})$.
\end{itemize}
\end{exam}

The example above is not an exception; in fact, in any nonreflexive Banach
space which is the dual of some Banach space one can construct a problem
$(\mathcal{P})$ having\ a nonzero gap with its biconjugate relaxation.

\begin{exam}
\label{examm}Assume\ that $X=Y^{\ast}$\ for some Banach space $Y.$ We define
the functions\ $f,$ $f_{y}:X\rightarrow\mathbb{R}$\ as
\[
f_{y}:=\left\langle \cdot,y\right\rangle ,\text{ }y\in B_{Y}\text{, and
}f:=\sup_{y\in B}f_{y}=\left\Vert \cdot\right\Vert _{Y^{\ast}}=\left\Vert
\cdot\right\Vert _{X}.
\]
We have
\[
f^{\ast\ast}=(\mathrm{I}_{B_{X^{\ast}}})^{\ast}=\mathrm{\sigma}_{B_{X^{\ast}}%
}=\left\Vert \cdot\right\Vert _{X^{\ast\ast}}\text{.}%
\]
At the same time,\ for every $y\in B_{Y}$ $(\subset B_{Y^{\ast\ast}%
}=B_{X^{\ast}})$ and $z\in X^{\ast\ast}$, by (\ref{moreau}) we have
\[
f_{y}^{\ast\ast}(z)=\liminf_{x\rightharpoonup z,\text{ }x\in X\text{ }}%
f_{y}(x)=\left\langle y,z\right\rangle ,
\]
and we deduce that
\[
\sup_{y\in B_{Y\text{ }}}f_{y}^{\ast\ast}=\sup_{y\in B_{Y}\text{ }%
}\left\langle y,\cdot\right\rangle =\mathrm{\sigma}_{B_{Y}}.
\]
We now consider the convex problems $(\mathcal{P}_{x^{\ast}}),$ $x^{\ast}\in
X^{\ast},$ given by
\[
(\mathcal{P}_{x^{\ast}})\text{ \ \ }%
\begin{array}
[c]{l}%
\inf~(u-\left\langle x^{\ast},x\right\rangle )\\
\text{s.t. }f_{y}(x)-u\leq0,\text{ }y\in B_{Y},\\
\text{ \ \ \ \ }(x,u)\in X\times\mathbb{R}.
\end{array}
\]
whose\ biconjugate relaxation is
\[
(\mathcal{P}_{x^{\ast}}^{\ast\ast})\text{ \ \ }%
\begin{array}
[c]{l}%
\inf~(u-\left\langle x^{\ast},z\right\rangle )\\
\text{s.t. }f_{y}^{\ast\ast}(z)-u\leq0,\text{ }y\in B_{Y},\\
\text{ \ \ \ \ }(z,u)\in X^{\ast\ast}\times\mathbb{R}.
\end{array}
\]
Observe that
\[
v(\mathcal{P}_{x^{\ast}})=\inf_{x\in X\text{ }}(f(x)-\left\langle x^{\ast
},x\right\rangle )=-f^{\ast}(x^{\ast}),
\]
while
\[
v(\mathcal{P}_{x^{\ast}}^{\ast\ast})=\inf_{x\in X\text{ }}(\sup_{y\in B_{Y}%
}f_{y}^{\ast\ast}(x)-\left\langle x^{\ast},x\right\rangle )=-(\sup_{y\in
B_{Y}}f_{y}^{\ast\ast})^{\ast}(x^{\ast}).
\]
If $v(\mathcal{P}_{x^{\ast}})=v(\mathcal{P}_{x^{\ast}}^{\ast\ast}),$ for all
$x^{\ast}\in X^{\ast},$ then we would have that $f^{\ast}(x^{\ast}%
)=(\sup_{y\in B_{Y}}f_{y}^{\ast\ast})^{\ast}(x^{\ast}),$ for all $x^{\ast}\in
X^{\ast}.$ Taking conjugates in the pair $(X^{\ast},X^{\ast\ast})$, we get
\[
f^{\ast\ast}(z)=\sup_{y\in B_{Y}}f_{y}^{\ast\ast}(z),\text{ for all }z\in
X^{\ast\ast};
\]
that is, $\left\Vert \cdot\right\Vert _{X^{\ast\ast}}=\mathrm{\sigma}_{B_{Y}%
}.$ Hence, 
\[
B_{Y}=B_{X^{\ast}}=B_{Y^{\ast\ast}},
\]
and we conclude that $Y^{\ast\ast}=Y$; in other words, $Y$ is reflexive. It follows that
$X=Y^{\ast}$ is reflexive too. We conclude that if $X$ is nonreflexive in this
example, then there must exist some $x_{0}^{\ast}\in X^{\ast}$ such that
\[
v(\mathcal{P}_{x_{0}^{\ast}}^{\ast\ast})<v(\mathcal{P}_{x_{0}^{\ast}}).
\]
Hence, a gap exists between $v(\mathcal{P}_{x_{0}^{\ast}}^{\ast\ast})$ and $v(\mathcal{P}_{x_{0}^{\ast}})$ despite the function $f$ is continuous and the Slater
condition holds.
\end{exam}

\section{Conclusions\label{sec5}}

\begin{enumerate}
\item The biconjugate relaxation $(\mathcal{P}^{\ast\ast})$ is adequate for
convex optimization problems with a finite number of constraints, since it
guarantees zero-duality gap under the Slater condition and the continuity of
the supremum of the constraint functions.

\item For infinite optimization problems, involving infinitely many
constraints, $(\mathcal{P}^{\ast\ast})$ may fail to satisfy a zero gap, as
shown in Section \ref{sec4}. In fact, a nontrivial example is given\ in
$c_{0}$ with a nonzero gap with its biconjugate relaxation. In addition, it is
shown that in any nonreflexive dual Banach space there exists a linear program
with infinitely many constraints presenting a gap with its biconjugate
relaxation. In this example, the supremum of the constraint functions is
continuous everywhere and the Slater condition is held.

\item To overcome this limitation, we propose a reinforced alternative to
$(\mathcal{P}^{\ast\ast});$ namely, $(\mathcal{P}_{\infty}^{\ast\ast}).$ We
prove that zero-duality gap holds between $(\mathcal{P})$ and $(\mathcal{P}%
_{\infty}^{\ast\ast}).$ Other variants related to\ $(\mathcal{P}_{\infty
}^{\ast\ast})$ are shown to be also useful for our purposes.

\item The assumptions underlying the theorems in this work are simple and
natural; namely, the Slater condition and the continuity of the supremum of
the constraint functions.

\item The analysis is carried out in the setting of\ infinite optimization
problems with a countable family of constraints. The case of general infinite
optimization problems with a possibly uncountable set of constraints can be
reduced to our setting when the underlying Banach space $X$ is separable.

\item By applying the biconjugate relaxation we establish zero-duality gap for
a new Fenchel-type dual $(\mathcal{D}^{\prime})$ of $(\mathcal{P}).$

\item The zero gap problem is closely related to the interchangeability of the
supremum and biconjugation\ operations, namely,
\begin{equation}
\left(  \sup_{t\in T}f_{t}\right)  ^{\ast\ast}=\sup_{t\in T}f_{t}^{\ast\ast},
\end{equation}
for a given family of convex functions $f_{t},$ $t\in T.$\ Such an identity
plays an important role in convex analysis (see, e.g., \cite{Ro12}), enabling
many subdifferential calculus of the supremum function $\sup_{t\in T}f_{t}$
(see, e.g., \cite{HLZ08}, \cite{LiNg11}; see also \cite{CHLBook} and the
references therein). The above relation has been recently investigated from a
subdifferential calculus perspective in \cite{CHL26-JOTA}, where it is shown
to be equivalent to strong subdifferential calculus rules for the supremum
function. Alternative characterizations based on closedness-type regularity
conditions have been developed in \cite{BG08}. The case of finite families has
been addressed\ in \cite{BW08, CHL26RACSAM, FS00, Za08}.
\end{enumerate}

\backmatter

%
%


\bmhead{Acknowledgements}

The research is supported by Grant PID2022-136399NB-C21
funded by MICIU/AEI/10.13039/501100011033 and by ERDF/EU, and Basal CMM
FB210005.

\section*{Declarations}


\begin{itemize}
\item Conflict of interest/Competing interests. Not applicable. 
\item Ethics approval and consent to participate. Not applicable.
\item Data availability. Not applicable. 
\end{itemize}

\end{document}